\documentclass[preprint,12pt]{elsarticle}

\usepackage{amssymb}

\usepackage{amsmath}                    
\usepackage{svg}                        
\usepackage{graphicx}                   
\usepackage{sistyle}                    
\usepackage[defaultlines=3,all]{nowidow}
\usepackage{booktabs}                   
\usepackage{tabularx}                   
\usepackage{adjustbox}                  
\usepackage{textcomp}                   
\usepackage{gensymb}                    
\usepackage{varwidth}                   
\usepackage{supertabular}               
\usepackage{multirow}                   
\usepackage{optidef} 				            
\usepackage{makecell}                   
\usepackage{listings}                   
\usepackage{tikz}                       
\usepackage{pgfplots}                   
\usepackage{mathrsfs}                   
\usepackage{tablefootnote}      
\usepackage{mdframed}
\usepackage{algorithm}
\usepackage{algpseudocode}
\usepackage{amsthm}                     
\usepackage{mathtools}                  

\usepackage{subcaption}

\usepgfplotslibrary{external}
\pgfplotsset{compat=newest}
\usetikzlibrary{positioning,calc,trees,angles,quotes}
\usetikzlibrary{shapes.geometric}
\usetikzlibrary{arrows.meta}
\usetikzlibrary{decorations.pathreplacing}
\usetikzlibrary{intersections}
\usetikzlibrary{fadings}
\usetikzlibrary{calc} 

\definecolor{MATLABblue}{RGB}{0,114,189} 
\definecolor{KITgreen}{RGB}{0,150,130}
\definecolor{mygreen}{RGB}{28,172,0}    
\definecolor{mylilas}{RGB}{170,55,241}

\journal{}

\begin{document}

  \begin{frontmatter}


    \title{A level set topology optimization theory based on Hamilton's principle} 

    \author[Oellerich]{Jan Oellerich}
    \author[Yamada]{Takayuki Yamada}

    \affiliation{organization={The University of Tokyo, Graduate School of Engineering, Institute of Engineering Innovation},
              addressline={Yayoi 2-11-16, Bunkyo-ku}, 
              city={Tokyo},
              postcode={113-8656}, 
              country={Japan}}

    \begin{abstract}
    In this paper, we present a novel framework for deriving the evolution equation of the level set function in topology optimization, departing from conventional Hamilton-Jacobi based formulations. The key idea is the introduction of an auxiliary domain, geometrically identical to the physical design domain, occupied by fictitious matter which is dynamically excited by the conditions prevailing in the design domain. By assigning kinetic and potential energy to this matter and interpreting the level set function as the generalized coordinate to describe its deformation, the governing equation of motion is determined via Hamilton's principle, yielding a modified wave equation. Appropriate combinations of model parameters enable the recovery of classical physical behaviors, including the standard and biharmonic wave equations. The evolution problem is formulated in weak form using variational methods and implemented in the software environment FreeFEM++. The influence of the numerical parameters is analyzed on the example of minimum mean compliance. The results demonstrate that topological complexity and strut design can be effectively controlled by the respective parameters. In addition, the method allows for the nucleation of new holes and eliminates the need for re-initializing the level set function. The inclusion of a damping term further enhances numerical stability. To showcase the versatility and robustness of our method, we also apply it to compliant mechanism design and a bi-objective optimization problem involving self-weight and compliance minimization under local stress constraints.
\end{abstract}

    \begin{graphicalabstract}
    \includegraphics{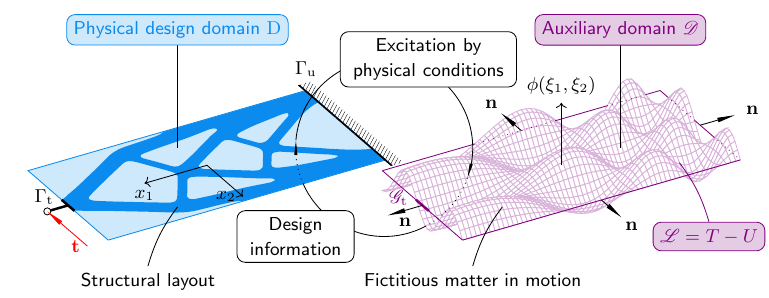}
\end{graphicalabstract}

    \begin{highlights}
      \item Unified framework for deriving the level set evolution in topology optimization.
      \item Evolution governed by Hamilton's principle, yielding a modified wave equation.
      \item Geometric complexity and strut shape can be controlled via numerical parameters.
      \item The method enables hole nucleation and does not require re-initialization.
      \item Robustness is demonstrated by varying the initial conditions.
    \end{highlights}

    \begin{keyword}
      Topology optimization, level set method, Hamilton's principle, wave equation, fictitious physical model
    \end{keyword}
  \end{frontmatter}

  \section{Introduction}
\label{sec:Introduction}
Topology optimization is an advanced method in the field of structural optimization and offers the highest degree of design freedom compared to so-called sizing optimization problems \citep{ALLAIRE2019}, in which the geometry of the component to be optimized is known in advance and parameterized accordingly, or shape optimization approaches \citep{HSU1994, SAITOU2005}, with which the shape of the component can be controlled. In topology optimization, the overall aim is to determine a preferred material distribution in a given domain that minimizes or maximizes a defined performance criterion while taking certain design constraints into account, depending on the respective application. Due to its flexibility the method it is not only limited to problems of structural mechanics but was also successfully demonstrated for electromagnetic \citep{CHOI2011}, thermal \citep{ZHUANG2007,WU2021,ONODERA2025-1,ONODERA2025-2} or gas flow problems \citep{GUAN2023,GUAN2024}.

Various methods for the treatment of topology optimization problems are reported in the literature while one important class covers density-based methods, in which a continuous material distribution is used to compute the optimal topology. Here, the Homogenization method proposed by Bends{\o}e and Kikuchi in \citep{BENDSOE1988} considers a periodic microstructural field with specific material properties while the density varies between 0 and 1 \citep{ALLAIRE2002}. Optimization is achieved by adjusting the microstructural parameters in order to obtain the desired material characteristics which makes the method preferable for designing materials possessing specific properties such as acoustic \citep{NOGUCHI2018, NOGUCHI2021} or electromagnetic metamaterials \citep{MURAI2023}, thermal related properties \citep{NAKAGAWA2023, GUO2024, OHENEBAAGYEKUM2025} or anisotropic composites \citep{KIM2020}. In the Solid Isotropic Material with Penalization method (SIMP), the design domain is divided into finite elements whose material stiffness can be controlled by a density variable and a penalization factor \citep{BENDSOE1999, SIGMUND2001, BENDSOE2004}. As a widely used method it has been applied to a broad range of optimization problems, such as multi-objective \citep{ZHANG2021} and multi-material problems \citep{NGUYEN2024, WAN2024} or manufacturing-related thickness control in component design \citep{ZUO2025}. 

In contrast to the previously outlined approaches, so-called boundary based methods employ front propagation of the structural boundary to create changes in the topology. Here, the level set method and the phase field method are important representatives while the latter introduces a smooth phase field function that is defined over the domain under consideration and determines the phase of each material point of the domain. This results in a diffuse interface between the single phases, while the thickness of the diffusive layer directly depends on the mesh refinement. In this context, the Van der Waals free energy of the system and the time-dependent evolution equation of the phase field function are related to each other while it is assumed that changes in the phase field function and the direction in which the free energy is minimized are in linear correspondence. This finally yields the Allen-Cahn equation which is employed to control the boundary movement \citep{TAKEZAWA2010}. While in traditional phase field methods both the objective functional and the interface energy are minimized, Takaki and Kato presented in \citep{TAKAKI2017} an approach that removes curvature effects due to the diffusion term and only focused on the minimization of the objective functional. The idea was to consider an alternative expression of the Laplacian yielding a one-dimensional phase field equation with respect to the interface normal direction. Seong et al. combined in \citep{SEONG2018} the reaction-diffusion equation with a modified conjugate gradient method (CGM) to increase the convergence rate ensuring numerical stability. Further work regarding the application of the reaction-diffusion equation in topology optimization can be taken from \citep{KIMC2020, JUNG2021}. A common drawback of the phase field method is that it is not capable to create new holes during the optimization process which is why it is frequently extended by additional techniques. To address this issue, Gao et al. combined in \citep{GAO2020} for instance the phase field method with the bidirectional evolutionary structural optimization (BESO) method which was utilized for the nucleation of new holes in the domain. An alternative approach to enable the creation of new holes was presented by Wang et al. in \citep{WANG2024} where the phase field method was combined with the optimality criteria approach originating from the field of machine learning which was employed for the removal of inefficient materials in the design domain. A combination of the phase field method and aspects of the density based SIMP method was presented by Jin et al. in \citep{JIN2024}. Here, the authors developed an adaptive algorithm that efficiently solves the benchmark problem of minimum mean compliance incorporating error estimators and proved its convergence. A similar approach which combines Filtered methods and the phase field method was conducted by Auricchio et al. in \citep{AURICCHIO2024}. In order to reduce computational costs regarding topology optimization in eigenfrequency problems, Hu et al. proposed in \citep{HU2023} a linearization approach for approximating eigenvalue problems in combination with the phase field method and demonstrated their approach on the example of a vibrating mechanical structures and photonic crystals. An extension of the phase field method to multi-scale structures incorporating an energy-based Homogenization method was investigated by Wang et al. in \citep{WANG2025}. Here, the Homogenization method was employed to assess the effectiveness of the microscopic structure with respect to the macroscopic performance of the structure while the phase field method was used to describe both the micro and macroscopic structure considering individual phase field functions. Comparative studies of different approaches containing investigations on parameters such as mesh sizes or diffusion coefficients were given by Kumar and Rakshit in \citep{KUMAR2024}. 

The level set method on the other hand is in the focus of our research and was first introduced by Osher and Sethian \citep{OSHER1988} to track the evolution of dynamically evolving interfaces and later transferred to the field of structural optimization \citep{SETHIAN2000}. The main idea is to implicitly represent the structural boundary using the isocontour of a higher-dimensional function, commonly referred to as the level set function. The implicit description of the component contour has the advantage that the component boundary is clearly represented and thus no diffuse transition zone is generated. The evolution of this function is governed by the Hamilton-Jacobi equation and is driven on the basis of sensitivity information obtained from the optimization process. Such an approach was proposed by Wang et al. in \citep{WANG2003} where the authors updated the level set function by means of an up-wind scheme being restricted by the Courant–Friedrichs–Lewy (CFL) conditions and required a re-initialization of the level set function after each iteration. Instead of using a finite difference scheme, Liu et al. employed a finite element solver based on the FEMLAB package for the reaction-diffusion equation to evolve the structural boundaries in minimum mean compliance problems \citep{LIU2005}. Here, the term containing the front propagation speed was considered the reaction term while an artificial dissipation term was added to the Hamilton-Jacobi equation. Jung et al. extended in \citep{JUNG2007} the level set method to the surface minimization of the zero level set with respect to triply-periodic surfaces that emerge in a wide range of systems including lipid-water systems or nanocomposites. The numerical implementation used the Hamilton-Jacobi equation and an auxiliary Newton iteration scheme to enforce the volume fraction constraint. A semi-implicit additive operator splitting scheme (AOS) was proposed by Luo et al. \citep{LUO2008} to numerically solve the Hamilton-Jacobi equation. The authors emphasized that the method is not restricted by the time step such as in explicit schemes due to the semi-implicit description, thus allowing an efficient implementation as it additionally refrains from a global initialization of the level set function. Zhou and Li \citep{ZHOU2008} applied the level set method to the topology optimization of steady state Navier-Stokes flow. In this context, the Hamilton-Jacobi equation was numerically solved by the implementation of two sets of meshes. To ensure that the level set function remains a signed distance function without altering its zero-level contour, the algorithm required frequent re-initialization. A similar work was proposed by Yamasaki et al. in \citep{YAMASAKI2010} who introduced a discretized signed distance function that was used as the design variables. Here, a re-initialization was also necessary to maintain the signed distance characteristic of the design variables. The authors applied their approach to the minimum mean compliance problem and compliant mechanism design, incorporating an additional perimeter constraint to address the typical issue of ill-posedness in structural optimization problems. A topology optimization approach that merges the level set method and elements of the phase field method was proposed by Yamada et al. in \citep{YAMADA2010}. Here, a regularization term acting as a fictitious interface energy was added to the Lagrangian functional that ensures smoothness of the level set function. An important feature was the possibility to control the resulting complexity of the topology by the regularization parameters as well as the capability of creating new holes during the optimization process thus combining both advantages of the phase field method and the level set method. Furthermore, the approach did not require a re-initialization of the level set function. Eventually, the method was successfully demonstrated on the examples of minimum mean compliance problem, compliant mechanism design and maximization of the lowest eigenfrequency. Zhou and Wang applied in \citep{ZHOU2012} a semi-Lagrangian level set method which implies the propagation of the level set function along a characteristic curve and allows for larger time steps compared to the up-wind scheme. The authors approximated the characteristic curve by using the first order Courant-Isaaacson-Rees formula to directly compute the level set function from the Hamilton-Jacobi equation and validated their method on the minimum mean compliance problem requiring a initially perforated design domain and re-initialization during the process. Xia et al. used in \citep{XIA2014} two level set functions to optimize both the structure and its support by means of forward Euler time differentiating to solve the  Hamilton-Jacobi equation. An alternative method to avoid a re-initialization of the level set function was given by Zhu et al. who introduced an additional energy functional that is added to the objective functional and used to obtain a generalized expression of the Hamilton-Jacobi equation \citep{ZHU2015}. A similar approach was reported by Jiang and Chen in \citep{JIANG2017}. A level set method based on the reaction-diffusion equation using a body fitted mesh was presented by Zhuang et al. in \citep{ZHUANG2021}. The proposed approach was implemented in Matlab and applied to benchmark problems such as minimum mean compliance and the design of compliant mechanisms, demonstrating its effectiveness. In their following work, the method was extended to nonlinear diffusion regularization \citep{ZHUANG2022}. Regarding further applications incorporating the reaction-diffusion equation in level set based topology optimization, reference is made to \citep{WANG2022, OKA2023, YU2024, TAJIMA2024, KUMAR2025, SOMA2025, YOSHIDA2025}.

A review of the literature on boundary-based topology optimization reveals that the underlying mathematical frameworks are largely confined to the reaction–diffusion equation, the Hamilton–Jacobi equation, or their combination. While numerous refinements have been proposed to enhance certain aspects such as preventing re-initialization or enabling hole nucleation these are often rooted in specific modifications rather than a unified theoretical basis. To address this limitation, the present work introduces a systematically derived evolution equation for the level set function, based on Hamilton's principle.

The remainder of this paper is organized as follows. Section~\ref{sec:mathematicalformulations} provides a brief overview of the baseline situation and commonly used formulations in level set based structural optimization. Building on this, an alternative perspective on the optimization problem is developed, which serves as the foundation for introducing an auxiliary domain and applying Hamilton's principle to derive the equation of motion for a fictitious matter. In this context, the level set function is interpreted as a generalized coordinate describing its deformation. Section~\ref{sec:numericalimplementation} presents the discretization of the evolution problem and its translation into a numerical update scheme based on its weak form. The proposed method is then applied to several design optimization problems, including an investigation of the influence of numerical parameters. Finally, Section~\ref{sec:conclusion} summarizes the main findings of this research.

  \section{Mathematical formulations}
\label{sec:mathematicalformulations}
\subsection{Baseline situation}
In structural optimization, we are interested in identifying that specific domain $\Omega^\prime\subset\mathbb{R}^N,N\in\{2,3\}$ containing solid material which minimizes a defined objective functional $J[\mathbf{u}(\Omega)]$ while the set of possible solutions is in general constrained by further design limitations and the governing equations of the physical system. In this context, consider first an isotropic and linear elastic body occupying an open material domain $\Omega$ being embedded in a design domain $\mathrm{D}\subset\mathbb{R}^N$ so that $\Omega\subseteq\mathrm{D}$. The boundary $\partial\Omega$ of the body is assumed to be smooth and defined as the union set of the boundary sections $\Gamma_\mathrm{u}$, $\Gamma_\mathrm{t}$, and the respective complement $\partial\Omega\setminus(\Gamma_\mathrm{u}\cup\Gamma_\mathrm{t})$. On $\Gamma_\mathrm{u}$ homogeneous Dirichlet boundary conditions are imposed while $\Gamma_\mathrm{t}$ is subject to external traction forces $\mathbf{t}$ and therefore leading to non-homogeneous Neumann boundary conditions. In contrast, the complement $\partial\Omega\setminus(\Gamma_\mathrm{u}\cup\Gamma_\mathrm{t})$ is assumed to be traction-free and hence subject to homogeneous Neumann boundary conditions. Summarizing these statements yields the set of partial differential equations
\begin{align}
    \left\{\begin{array}{rll}
        -\mathrm{div}(\mathbf{C}:\varepsilon(\mathbf{u})) &= \mathbf{b} & \text{in}\quad\Omega \\
        \mathbf{u} &= \mathbf{0} & \text{on}\quad\Gamma_{\mathrm{u}} \\
        (\mathbf{C}:\varepsilon(\mathbf{u}))\mathbf{n} &= \mathbf{t} & \text{on}\quad\Gamma_{\mathrm{t}} \\
        (\mathbf{C}:\varepsilon(\mathbf{u}))\mathbf{n} &= \mathbf{0} & \text{on}\quad\partial\Omega\setminus(\Gamma_\mathrm{u}\cup\Gamma_\mathrm{t})
        \end{array}\right.
        \label{eq:governingequations}
\end{align}
describing the equilibrium condition of the body in its strong formulation considering body forces $\mathbf{b}$. Here, $\mathbf{C} = 2\mu_\mathrm{c}\mathbb{I}+\lambda_\mathrm{c}\mathbf{I}\otimes\mathbf{I}$ ($\mathbb{I}$ and $\mathbf{I}$ as fourth and second rank identity tensors) denotes the elasticity tensor with Lam\'e constants $\mu_\mathrm{c} = E / (2(1+\nu))$ and $\lambda_\mathrm{c} = E\nu/((1+\nu)(1-2\nu))$ in plane strain assumption taking Young's modulus $E$ and Possion's ratio $\nu$ into account. Furthermore, $\varepsilon(\mathbf{u}) = (\nabla\mathbf{u} + (\nabla\mathbf{u})^\top)/2$ relates to the linearized strain tensor being applied to the displacement field $\mathbf{u}(\Omega)$ which represents the state variable. Using Eq.~(\ref{eq:governingequations}) and considering an additional inequality constraint functional $G[\mathbf{u}]$, the original structural optimization problem is given as
\begin{customopti}|s|
    {inf}{\Omega\subseteq\mathrm{D}\subset\mathbb{R}^N}{
    J[\mathbf{u}] = \int_\Omega j(\mathbf{u})\,\mathrm{d}\Omega}{}{}{}
    \addConstraint{G[\mathbf{u}]}{\leq 0}{\quad\text{in}\quad\Omega}{}
    \addConstraint{-\mathrm{div}(\mathbf{C}:\varepsilon(\mathbf{u}))}{=\mathbf{b}}{\quad\text{in}\quad\Omega}{}
    \addConstraint{\mathbf{u}}{=\mathbf{0}}{\quad\text{on}\quad\Gamma_\mathrm{u}}{}
    \addConstraint{(\mathbf{C}:\varepsilon(\mathbf{u}))\mathbf{n}}{=\mathbf{t}}{\quad\text{on}\quad\Gamma_\mathrm{t}}{}
    \addConstraint{(\mathbf{C}:\varepsilon(\mathbf{u}))\mathbf{n}}{=\mathbf{0}}{\quad\text{on}\quad\partial\Omega\setminus(\Gamma_\mathrm{u}\cup\Gamma_\mathrm{t})}{}.
    \label{eq:originalproblem}
\end{customopti}
The formulation in Eq.~(\ref{eq:originalproblem}) clearly defines the optimization problem, however, the solution turns out to be challenging since an optimum material domain $\Omega^\prime$ must be determined directly. Therefore, the problem is transformed in the next step into a material distribution problem by introducing the characteristic function
\begin{align*}
    \chi_\Omega(\mathbf{x}) = 
    \left\{\begin{array}{rl}
        1, & \mathbf{x}\in\Omega\cup\partial\Omega\\[10 pt]
        0, & \mathbf{x}\in\mathrm{D}\setminus(\Omega\cup\partial\Omega)
    \end{array}\right.
\end{align*}
which is an element of the bounded Lebesgue space $L^\infty(\Omega)$ and enables the formal description of any topology. This allows the optimization problem to be extended to the design domain $\mathrm{D}$ and which then follows to
\begin{customopti}|s|
    {inf}{\chi_\Omega\in L^\infty(\Omega)}{
    J[\mathbf{u},\chi_\Omega] = \int_\mathrm{D}j(\mathbf{u})\chi_\Omega\,\mathrm{d}\Omega}{}{}{}
    \addConstraint{G[\mathbf{u},\chi_\Omega]}{\leq 0}{\quad\text{in}\quad\mathrm{D}}{}
    \addConstraint{-\mathrm{div}(\mathbf{C}:\varepsilon(\mathbf{u}))}{=\mathbf{b}}{\quad\text{in}\quad\Omega}{}
    \addConstraint{\mathbf{u}}{=\mathbf{0}}{\quad\text{on}\quad\Gamma_\mathrm{u}}{}
    \addConstraint{(\mathbf{C}:\varepsilon(\mathbf{u}))\mathbf{n}}{=\mathbf{t}}{\quad\text{on}\quad\Gamma_\mathrm{t}}{}
    \addConstraint{(\mathbf{C}:\varepsilon(\mathbf{u}))\mathbf{n}}{=\mathbf{0}}{\quad\text{on}\quad\partial\Omega\setminus(\Gamma_\mathrm{u}\cup\Gamma_\mathrm{t})}{}.
    \label{eq:characteristicfunctionproblem}
\end{customopti}
Although the characteristic function $\chi_\Omega$ ensures integrability, its implementation is well known to cause ill-posedness \citep{ALLAIRE2002}, as solutions may exhibit discontinuities over infinitesimal regions. To address this issue, various approaches have been proposed in the literature, including the previously outlined Homogenization method. Another important approach is the level set method, which introduces a higher-dimensional function $\phi\in H^\alpha(\mathrm{D})$ as its key concept. In this context, $H^\alpha(\mathrm{D})$ with $\alpha\geq 1$ denotes the Sobolev space over the domain $\mathrm{D}$ containing weak differentiable functions of order $\alpha$ that are square-integrable. Its norm is induced by its inner product. The level set function $\phi:\mathbb{R}^N\rightarrow[-1,1]$ is defined as 
\begin{align}
    \left\{\begin{array}{rl}
        1 \geq\phi(\mathbf{x}) > 0, & \mathbf{x}\in\Omega\\
        \phi(\mathbf{x}) = 0, & \mathbf{x}\in\partial\Omega \\
        -1\leq\phi(\mathbf{x}) < 0, & \mathbf{x}\in\mathrm{D}\setminus(\Omega\cup\partial\Omega)
    \end{array}\right.
    \label{eq:levelsetfunctiondefinition}
\end{align}
and implicitly describes the structural boundary through its zero isocontour. Thus, the structural optimization problem can be reformulated as follows:
\begin{customopti}|s|
    {inf}{\phi\in H^\alpha(\mathrm{D})}{
    J[\mathbf{u},\chi(\phi)] = \int_\mathrm{D}j(\mathbf{u})\chi(\phi)\,\mathrm{d}\Omega}{}{}{}
    \addConstraint{G[\mathbf{u},\chi(\phi)]}{\leq 0}{\quad\text{in}\quad\mathrm{D}}{}
    \addConstraint{-\mathrm{div}(\mathbf{C}:\varepsilon(\mathbf{u}))}{=\mathbf{b}}{\quad\text{in}\quad\Omega}{}
    \addConstraint{\mathbf{u}}{=\mathbf{0}}{\quad\text{on}\quad\Gamma_\mathrm{u}}{}
    \addConstraint{(\mathbf{C}:\varepsilon(\mathbf{u}))\mathbf{n}}{=\mathbf{t}}{\quad\text{on}\quad\Gamma_\mathrm{t}}{}
    \addConstraint{(\mathbf{C}:\varepsilon(\mathbf{u}))\mathbf{n}}{=\mathbf{0}}{\quad\text{on}\quad\partial\Omega\setminus(\Gamma_\mathrm{u}\cup\Gamma_\mathrm{t})}{}.
    \label{eq:levelsetfunctionproblem}
\end{customopti}
Furthermore, the weak form of the equilibrium condition is derived by means of variational methods. For this purpose, we first define the function space
\begin{align*}
    \mathscr{U} = \{\mathbf{u} = u_i\mathbf{e}_i, u_i\in H^1(\mathrm{D})\mid\mathbf{u} = \mathbf{0}\quad\text{on}\quad\Gamma_\mathrm{u}\}
\end{align*}
and multiply Eq.~(\ref{eq:governingequations}) by a suitable test function $\mathbf{v}\in\mathscr{U}$. Using the divergence theorem yields
\begin{align*}
    R[\mathbf{u},\chi(\phi);\mathbf{v}] = a(\mathbf{u},\mathbf{v},\chi(\phi)) - l(\mathbf{v},\chi(\phi)) = 0
\end{align*}
with linear operators
\begin{align*}
    a(\mathbf{u},\mathbf{v},\chi(\phi)) &= \int_\mathrm{D}\varepsilon(\mathbf{u}):\mathbf{C}:\varepsilon(\mathbf{v})\chi(\phi)\,\mathrm{d}\Omega,\\
    l(\mathbf{v},\chi(\phi)) &= \int_\mathrm{D}\mathbf{b}\cdot\mathbf{v}\chi(\phi)\,\mathrm{d}\Omega + \int_{\Gamma_\mathrm{t}}\mathbf{t}\cdot\mathbf{v}\,\mathrm{d}\Gamma
\end{align*}
in level set formulation. Note that $R[\mathbf{u},\chi(\phi);\mathbf{v}]$ represents the weak form of the equilibrium equation, which must vanish for a valid solution and may vary depending on the specific problem settings. Furthermore, it is convenient to express the characteristic function by means of a Heaviside function defined as
\begin{align}
    \Theta(\phi) \vcentcolon = \left\{\begin{array}{rl}
        1, &\quad\text{if}\quad\phi(\mathbf{x}) \geq 0\\[10 pt]
        0, &\quad\text{otherwise}
    \end{array}\right.
    \label{eq:heavisidefunction}
\end{align}
which gives the structural optimization problem
\begin{customopti}|s|
    {inf}{\phi\in H^\alpha(\mathrm{D})}{
    J[\mathbf{u},\Theta] = \int_\mathrm{D}j(\mathbf{u})\Theta\,\mathrm{d}\Omega}{}{}{}
    \addConstraint{G[\mathbf{u},\Theta]}{\leq 0}{\quad\text{in}\quad\mathrm{D}}{}
    \addConstraint{R[\mathbf{u},\Theta;\mathbf{v}]}{=0}{\quad\text{in}\quad\mathrm{D}}{}.
    \label{eq:baselineproblem}
\end{customopti}
Following the conventional procedure, the constrained optimization problem from Eq.~(\ref{eq:baselineproblem}) is usually transformed into an unconstrained problem using the Lagrangian method by adding the constraint functionals to the actual objective functional and assigning an additional multiplier $\lambda$ to the inequality constraint resulting in the Lagrangian functional
\begin{align}
    \mathcal{L}[\mathbf{u},\Theta;\mathbf{v},\lambda] = J[\mathbf{u},\Theta] + \lambda G[\mathbf{u},\Theta] - R[\mathbf{u},\Theta;\mathbf{v}].
    \label{eq:lagrangian}
\end{align} 
A function $\phi^\prime\in H^\alpha(\mathrm{D})$ which satisfies the necessary criteria of optimality known as the KKT conditions
\begin{align*}
    \left\langle\frac{\mathrm{d}\mathcal{L}}{\mathrm{d}\phi},\delta\phi\right\rangle = 0, && \lambda G[\mathbf{u},\chi(\phi)] = 0, && \lambda \geq 0, && G[\mathbf{u},\chi(\phi)] \leq 0
\end{align*}
is considered an optimal solution candidate. Here, the first term represents the Frech\'et derivative of the Lagrangian functional with respect to $\phi$ in the direction $\delta\phi$. A direct computation of $\phi^\prime$ yielding an optimum solution is non-trivial; therefore the solution is approached iteratively. For this purpose, a fictitious time $t$ is introduced and it is postulated that the level set function depends implicitly on this time, i.e. $\phi = \phi(\mathbf{x}, t)$. To track the evolution of the structural boundary, the level set function is differentiated with respect to $t$, which leads to the Hamilton-Jacobi equation
\begin{align}
    \frac{\partial\phi(\mathbf{x},t)}{\partial t} + \frac{\partial\mathbf{x}}{\partial t}|\nabla\phi(\mathbf{x},t)| = 0. 
    \label{eq:hamiltonjacobiequation}
\end{align}
By solving Eq.~(\ref{eq:hamiltonjacobiequation}) for each time step, the structural boundary is iteratively updated.
\subsection{Derivation of equation of motion}
The Hamilton-Jacobi equation serves as the foundation for topology optimization in current level set-based methods. As shown in \citep{YAMADA2010}, it can be transformed into a reaction-diffusion equation through appropriate modifications, thereby defining the mathematical characteristics of the evolution of the structural boundary. Consequently, the possible outcomes of topology optimization are inherently dictated by this mathematical formulation, raising the question of whether this description can be further extended. 

\begin{figure}[ht]
    \centering
    \includegraphics[width=1.\textwidth]{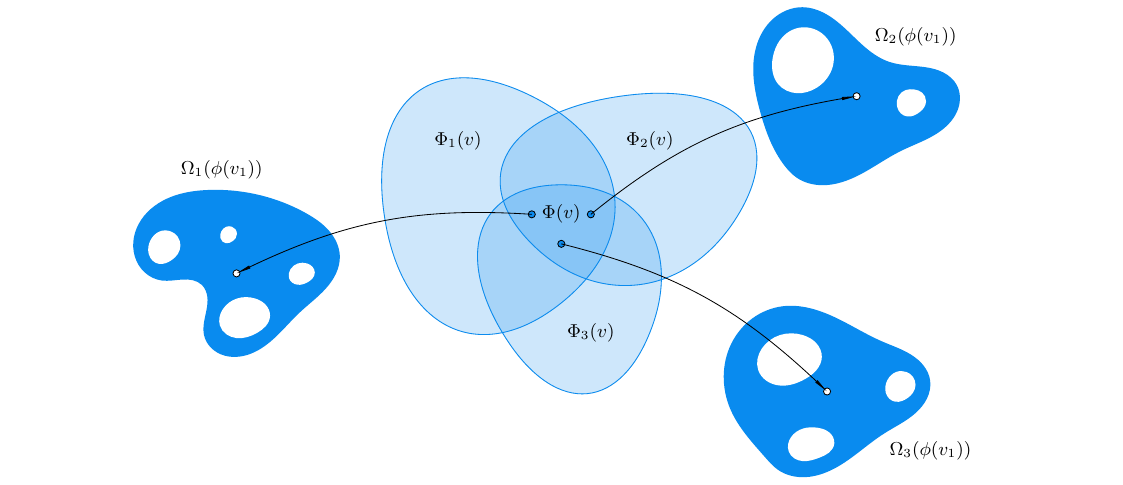}
    \caption{Construction of the set $\Phi(v_1)$ containing admissible solutions $\phi(v_1)$ leading to different structural configurations $\Omega_i(\phi(v_1))$.}
    \label{fig:admissiblesolutions}
\end{figure}

\begin{figure}[ht]
    \centering
    \includegraphics[width=1.\textwidth]{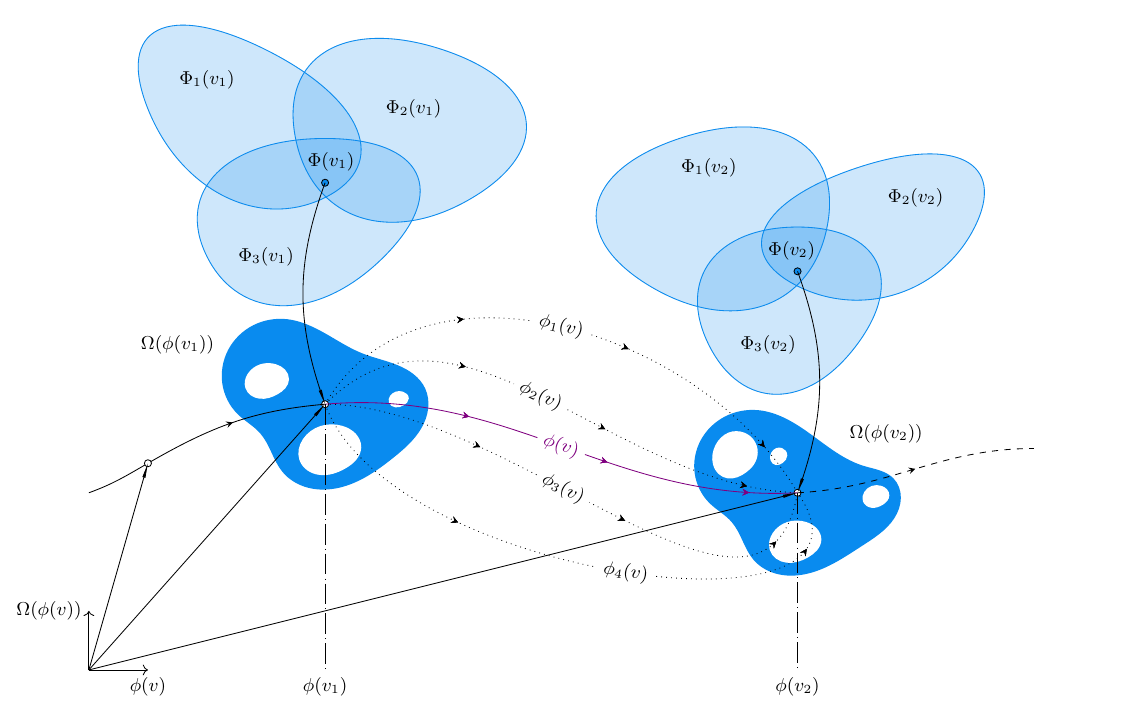}
    \caption{Evolution of the structural layout described by admissible trajectories $\phi_i(v)$ within the interval $[v_1,v_2]$.}
    \label{fig:evolution}
\end{figure}

\begin{figure}[ht]
    \centering
    \includegraphics[width=1.0\textwidth]{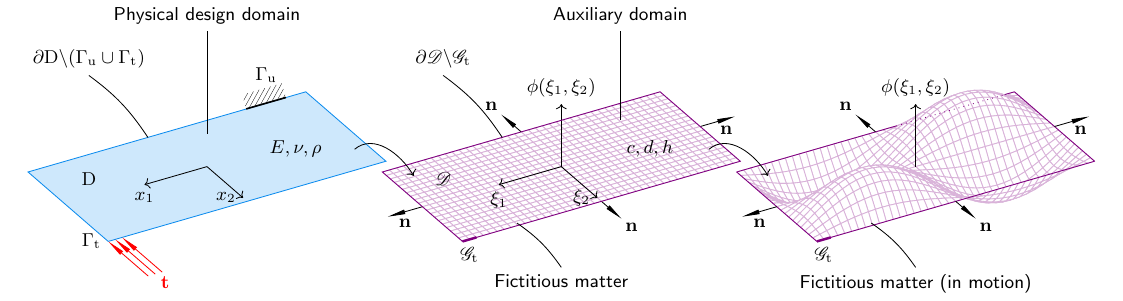}
    \caption{Auxiliary domain being occupied with fictitious matter.}
    \label{fig:auxiliarydomain}
\end{figure}

\begin{figure}[ht]
    \centering
    \includegraphics[width=1.\textwidth]{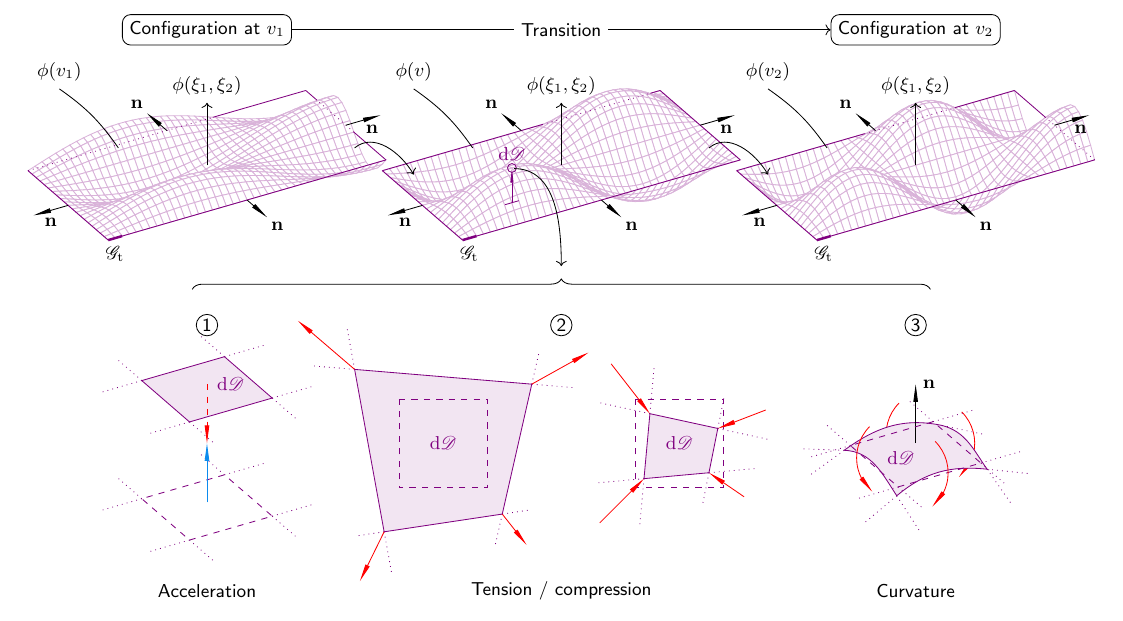}
    \caption{Transition of the fictitious matter within the interval $[v_1,v_2]$.}
    \label{fig:transition}
\end{figure}

For this purpose, we deliberately take an alternative approach and assume in a first instance that an admissible solution to the optimization problem defined in Eq.~(\ref{eq:baselineproblem}) exists. Using this hypothesis, we can deduce that the level set function apparently yields a certain value $v_1$ of the objective functional. Thus, the level set function must be a solution to the following equations: 
\begin{align*}
    G[\mathbf{u},\Theta]\leq 0, && R[\mathbf{u},\Theta;\mathbf{v}] = 0, && J[\mathbf{u},\Theta] - v_1 = 0.
\end{align*}
Next, we assume that more than one level set function exist which fulfill the inequality constraint as well as the equilibrium condition and additionally leads to the identical functional value $v_1$. With this in mind we are able to construct the set of admissible solutions $\Phi$ depending on the functional value $v_1$ as 
\begin{align*}
    \Phi(v_1) = \Phi_1(v_1)\cap\Phi_2(v_1)\cap\Phi_3(v_1) 
\end{align*}
considering the single sets
\begin{align*}
    \Phi_1 &\vcentcolon = \{\phi\in H^\alpha(\mathrm{D})\mid G[\mathbf{u},\Theta]\leq 0\}, \\
    \Phi_2 &\vcentcolon = \{\phi\in H^\alpha(\mathrm{D})\mid R[\mathbf{u},\Theta;\mathbf{v}] = 0\}, \\
    \Phi_3 &\vcentcolon = \{\phi\in H^\alpha(\mathrm{D})\mid J[\mathbf{u},\Theta] - v_1 = 0\}.
\end{align*}
Following this line of thought, we find as a consequence that the admissible solutions contained in $\Phi(v_1)$ may also differ in their resulting topology $\Omega(\phi(v_1))$ as qualitatively illustrated in Fig.~(\ref{fig:admissiblesolutions}). 
Since there are obviously several solutions to choose from, the fundamental question now arises as to how select a suitable solution from the solution set. To answer this, we propose the approach of varying the required value $v_1$ and consider a neighboring condition which yields the objective functional value $v_2$, see Fig.~(\ref{fig:evolution}).
In this context, we observe that the set of admissible solutions, i.e., $\Phi(v_2)$, changes as well, thereby altering the space of possible level set functions. Instead of identifying individual solutions, we now focus on the trajectory an admissible solution takes from $v_1$ to $v_2$ and postulate that this transition occurs continuously. Furthermore, we assume a natural behavior of the associated level set function as it undergoes that transition. Recall that the level set function is an auxiliary construct defined over the physical domain, employed to implicitly describe the structural boundary. This enables us to decouple the level set function from the physical domain and instead consider its evolution within an auxiliary domain $\mathscr{D}$, which is geometrically identical to the physical domain $D$. In the subsequent step, we seek a systematic description of the level set function's evolution. To this end, we propose that the auxiliary domain is occupied by fictitious matter that is not at rest but in motion, with the level set function serving as a state variable that characterizes the deflection of this fictitious matter, see Fig.~(\ref{fig:auxiliarydomain}).
By evaluating the zero isocontour of the level set function, we obtain information regarding the structural boundary and thus the resulting topology. This enables us to use the motion of the fictitious matter to describe the transition of the structural layout determined by $\phi(v_1)$ to $\phi(v_2)$ by considering its configurations in both states. As illustrated in Fig.~(\ref{fig:evolution}), multiple possible paths $\phi_i(v)$ (indicated by dotted lines) may exist along which the level set function could evolve. However, under the assumption of a natural behavior, we postulate that a specific trajectory $\phi(v)$ is preferred among all admissible candidates. This assumption allows us to employ the fundamental principle of least action from theoretical physics, which states that physical systems evolve in such a way that the associated action attains an extremum. In this context, Hamilton's principle is a specific formulation of the least action principle, asserting that the time integral of the Lagrangian function $\mathscr{L} = T - U$ takes an extreme value, where $T$ denotes the kinetic energy and $U$ the potential energy of the system \citep{LANCZOS1986,FOX1987,BEDFORD2021}. By transferring Hamilton's principle to our setting, we treat the level set function as the generalized coordinate and consider the evolution from $v_1$ to $v_2$ via the action integral
\begin{align*}
    \mathcal{S}[\phi(v)] = \int_{v_1}^{v_2}\mathscr{L}(\phi(v),\partial\phi(v)/\partial v,v)\,\mathrm{d}v.
\end{align*}
Let $\phi(v)$ denote the true trajectory between the two states $\phi(v_1)$ and $\phi(v_2)$ defined over the interval $[v_1,v_2]$ and let $\delta\phi(v)$ be an admissible variation that vanishes at the endpoints of the interval, i.e., $\delta(v_1)=\delta\phi(v_2)=0$. In addition, let $\epsilon$ be a small constant, then we derive the first variation as follows:
\begin{align}
    \notag\left\langle\frac{\partial\mathcal{S}}{\partial\phi},\delta\phi\right\rangle &= \frac{\mathrm{d}}{\mathrm{d}\epsilon}\left[\int_{v_1}^{v_2}\mathscr{L}(\phi + \epsilon\delta\phi,\partial(\phi + \epsilon\delta\phi)/\partial v, v)\,\mathrm{d}v\right]_{\epsilon=0}\\
    &= \int_{v_1}^{v_2}\left(\frac{\partial\mathscr{L}}{\partial\phi}\delta\phi + \frac{\partial\mathscr{L}}{\partial(\partial\phi/\partial v)}\frac{\partial\delta\phi}{\partial v}\right)\,\mathrm{d}v
    \label{eq:firstvariation}
\end{align}
For the second term in Eq.~(\ref{eq:firstvariation}) we perform integration by parts and obtain
\begin{align*}
    \int_{v_1}^{v_2}\frac{\partial\mathscr{L}}{\partial(\partial\phi/\partial v)}\frac{\partial\delta\phi}{\partial v}\,\mathrm{d}v = \left[\frac{\partial\mathscr{L}}{\partial(\partial\phi/\partial v)}\frac{\partial\delta\phi}{\partial v}\right]_{v_1}^{v_2} - \int_{v_1}^{v_2}\frac{\partial}{\partial v}\left(\frac{\partial\mathscr{L}}{\partial(\partial\phi/\partial v)}\right)\delta\phi\,\mathrm{d}v.
\end{align*}
Here, the first terms vanishes, so that we finally obtain
\begin{align}
    \left\langle\frac{\partial\mathcal{S}}{\partial\phi},\delta\phi\right\rangle = \int_{v_1}^{v_2}\left\{\frac{\partial\mathscr{L}}{\partial\phi} - \frac{\partial}{\partial v}\left(\frac{\partial\mathscr{L}}{\partial(\partial\phi/\partial v)}\right)\right\}\delta\phi\,\mathrm{d}v.
\end{align}
Since the action assumes an extremum, its first variation must vanish, i.e., $\delta\mathcal{S}=0$. This condition is satisfied if and only if the expression in the curly brackets vanishes, leading to the equation
\begin{align}
    \frac{\partial}{\partial v}\left(\frac{\partial\mathscr{L}}{\partial(\partial\phi/\partial v)}\right) - \frac{\partial\mathscr{L}}{\partial\phi} \overset{!}{=} 0
    \label{eq:eulerlagrange}
\end{align}
which is commonly known as the Euler Lagrangian equation. It is worth emphasizing that, unlike standard level set approaches where a fictitious time parameter $t$ governs the evolution, our method employs the parameter $v$, thereby describing the evolution in terms of changes in the objective functional. This effectively characterizes the motion as a trajectory along decreasing functional values. The resulting equation describes the free motion of the fictitious matter. However, we additionally assume the presence of a dissipative term denoted by $Q^\ast$ to ensure that the fictitious matter eventually reaches a stable condition. Therefore, we extend Eq.~(\ref{eq:eulerlagrange}) by $Q^\ast$ and obtain the new form
\begin{align}
    \frac{\partial}{\partial v}\left(\frac{\partial\mathscr{L}}{\partial(\partial\phi/\partial v)}\right) - \frac{\partial\mathscr{L}}{\partial\phi} = Q^\ast
    \label{eq:eulerlagrangedissipative}
\end{align}
which will be the foundation for the further course of the work. As stated earlier, the Lagrangian function $\mathscr{L}$ consists of the kinetic and potential energy of the fictitious matter. To identify appropriate expressions for these terms, we theoretically analyze the fictitious matter in its deformed state, as illustrated in Fig.~(\ref{fig:transition}). 
We consider an infinitesimal surface element $\mathrm{d}\mathscr{D}$ and assume that it experiences an acceleration during the transition. From this, we derive the kinetic energy of the entire fictituous matter which is proportional to the square of the rate of change with respect to $v$, given by
\begin{align*}
    T[\partial\phi/\partial v] = \frac{1}{2}\int_\mathscr{D}\left(\frac{\partial\phi}{\partial v}\right)^2\,\mathrm{d}\mathscr{D}
\end{align*}
where level set function serves to describe the deflection. For the potential energy we consider finite deformations of the fictitious matter and focus on a surface element $\mathrm{d}\mathscr{D}$. First, we note that the surface element undergoes compression and tension. Thus, the corresponding potential energy is proportional to the squared absolute value of the gradient of the level set function $\phi$, multiplied by a constant $c$ that represents the resistance to compression and tension. Since we also assume finite deformations we further observe that the surface element experiences curvature which is proportional to the second spatial derivative of $\phi$, multiplied by a factor $h^2$ that represent the resistance to curvature. Thus, we obtain 
\begin{align*}
    U[\phi] = U_1[\phi] + U_2[\phi] = \frac{1}{2}\int_\mathscr{D}(c|\nabla\phi|)^2\,\mathrm{d}\mathscr{D} + \frac{1}{2}\int_\mathscr{D}(h^2\nabla^2\phi)^2\,\mathrm{d}\mathscr{D}
\end{align*}
for the total potential energy. In the following, the single components of Eq.~(\ref{eq:eulerlagrangedissipative}) based on the expressions for the kinetic and potential energy are to be derived. The partial functional derivative of $\mathscr{L}$ with respect to $\partial\phi/\partial v$ follows to
\begin{align*}
    \left\langle\frac{\partial\mathscr{L}}{\partial(\partial\phi/\partial v)},\delta\left(\frac{\partial\phi}{\partial v}\right)\right\rangle &= \frac{\mathrm{d}}{\mathrm{d}\epsilon}\left[\frac{1}{2}\int_\mathscr{D}\left(\frac{\partial\phi}{\partial v} + \epsilon\delta\left(\frac{\partial\phi}{\partial v}\right)\right)^2\,\mathrm{d}\mathscr{D}\right]_{\epsilon=0} \\
    &= \int_\mathscr{D}\frac{\partial\phi}{\partial v}\delta\left(\frac{\partial\phi}{\partial v}\right)\,\mathrm{d}\mathscr{D}
\end{align*}
which first yields the derivative
\begin{align*}
    \frac{\partial\mathscr{L}}{\partial(\partial\phi/\partial v)} = \frac{\partial\phi}{\partial v}.    
\end{align*}
By additional differentiating with respect to $v$, we obtain the first term
\begin{align}
    \frac{\partial}{\partial v}\left(\frac{\partial\mathscr{L}}{\partial(\partial\phi/\partial v)}\right) = \frac{\partial^2\phi}{\partial v^2}.
\end{align}
Next, the partial derivative of the Lagrangian functional $\mathscr{L}$ with respect to $\phi$ is to be determined. Since the kinetic energy does not depend on $\phi$, we concentrate on the potential energy considering the boundary conditions
\begin{align*}
    \left\{\begin{array}{rl}
         \phi &= 1  \quad\text{on}\quad \mathscr{G}_\mathrm{t}\\[10 pt]
         \nabla\phi\cdot\mathbf{n} &=0 \quad\text{on}\quad \partial\mathscr{D}\setminus\mathscr{G}_\mathrm{t} 
    \end{array}\right..
\end{align*}
Here, we obtain for the first term
\begin{align*}
    \left\langle\frac{\partial U_1}{\partial \phi},\delta\phi\right\rangle &= \frac{\mathrm{d}}{\mathrm{d}\epsilon}\left[\frac{1}{2}\int_\mathscr{D}c^2\nabla(\phi+\epsilon\delta\phi)\cdot\nabla(\phi+\epsilon\delta\phi)\,\mathrm{d}\mathscr{D}\right]_{\epsilon=0} \\
    &= -\int_\mathscr{D}c^2\nabla^2\phi\delta\phi\,\mathrm{d}\mathscr{D}
\end{align*}
where we identify the partial derivative accordingly as
\begin{align*}
    \frac{\partial U_1}{\partial\phi} = - c^2\nabla^2\phi.
\end{align*}
The second term is derived in an analogous way where we conveniently use the substitution $\psi=\nabla^2\phi$. Then, the calculation proceeds as outlined above, thus we obtain after re-substitution the partial derivative as follows:
\begin{align*}
    \frac{\partial U_2}{\partial\phi} = - h^4\nabla^4\phi.
\end{align*}
In the next step, we determine the term $Q^\ast$ by introducing a dissipative functional $\mathscr{R}$, defined as 
\begin{align*}
    \mathscr{R}[\phi,\partial\phi/\partial v] = \mathscr{R}_1[\phi] + \mathscr{R}_2[\partial\phi/\partial v]
\end{align*}
where
\begin{align*}
    \mathscr{R}_1[\phi] &= G^\ast[\mathbf{u},\Theta;\lambda] + (J[\mathbf{u},\Theta] - v) - R[\mathbf{u},\Theta;\mathbf{v}], \\
    \mathscr{R}_2[\partial\phi/\partial v] &= \frac{1}{2}\int_\mathscr{D}d\left(\frac{\partial\phi}{\partial v}\right)^2\,\mathrm{d}\mathscr{D}.
\end{align*}
This formulation is inspired by the structure of a Rayleigh dissipation function \citep{LEMOS1991, ROY2007} and accounts for both dissipative effects and external influences, which are governed by the conditions in the design domain. These cover in particular the information obtained from the constraint and objective functionals, as well as the equilibrium state in $\mathscr{R}_1[\phi]$. Here, the constraint functional is to be treated via the Augmented Lagrangian method. The second term $\mathscr{R}_2[\partial\phi/\partial v]$ corresponds to a kinetic energy related dissipation multiplied by a factor $d$. At this point, it is worthy to note that the dissipative functional serves as the link between the design and auxiliary domain. In the following, we postulate the relation
\begin{align*}
    Q^\ast = -\frac{\partial\mathscr{R}}{\partial\phi} - \frac{\partial\mathscr{R}}{\partial(\partial\phi/\partial v)}
\end{align*}
and calculate the partial derivatives accordingly while the first term results to
\begin{align*}
    \frac{\partial\mathscr{R}}{\partial\phi} &= \left[\frac{\partial G^\ast}{\partial\mathbf{u}}\frac{\partial\mathbf{u}}{\partial\Theta}\delta(\phi) + \frac{\partial G^\ast}{\partial\Theta}\delta(\phi)\right] + \left[\frac{\partial J}{\partial\mathbf{u}}\frac{\partial\mathbf{u}}{\partial\Theta}\delta(\phi)+\frac{\partial J}{\partial\Theta}\delta(\phi)\right] - \dots\\
    &\dots-\left[\frac{\partial R}{\partial\mathbf{u}}\frac{\partial\mathbf{u}}{\partial\Theta}\delta(\phi) + \frac{\partial R}{\partial\Theta}\delta(\phi)\right]
\end{align*}
In this context, $\delta(\phi)$ refers to the Delta Dirac function. Next, we separate the implicit and explicit terms \citep{MICHALERIS1994} and obtain
\begin{align*}
    \frac{\partial\mathscr{R}}{\partial\phi} = \left[\frac{\partial G^\ast}{\partial\mathbf{u}} + \frac{\partial J}{\partial\mathbf{u}} - \frac{\partial R}{\partial\mathbf{u}}\right]\frac{\partial\mathbf{u}}{\partial\Theta}\delta(\phi) + \left[\frac{\partial G^\ast}{\partial\Theta} + \frac{\partial J}{\partial\Theta} - \frac{\partial R}{\partial\Theta}\right]\delta(\phi).    
\end{align*}
From this, we identify the adjoint equation which must vanish so that the implicit terms $\partial\mathbf{u}/\partial\Theta$ vanish. In its weak formulation, it is then given as
\begin{align}
    \mathcal{A}:\quad\left\langle\frac{\partial G^\ast}{\partial\mathbf{u}},\delta\mathbf{u}\right\rangle + \left\langle\frac{\partial J}{\partial\mathbf{u}},\delta\mathbf{u}\right\rangle = \left\langle\frac{\partial R}{\partial\mathbf{u}},\delta\mathbf{u}\right\rangle
    \label{eq:adjointequation}
\end{align}
considering an appropriate test function $\delta\mathbf{u}\in\mathscr{U}$. Let us now assume, that this condition is fulfilled, then the second term reduces to the form
\begin{align*}
    \frac{\partial\mathscr{R}}{\partial\phi} = \left[\frac{\partial G^\ast}{\partial\Theta} + \frac{\partial J}{\partial\Theta} - \frac{\partial R}{\partial\Theta}\right]\delta(\phi)
\end{align*}
where the terms in brackets are conveniently summarized to
\begin{align}
    f(\mathbf{u},\Theta;\mathbf{v},\lambda) = \frac{\partial G^\ast}{\partial\Theta} + \frac{\partial J}{\partial\Theta} - \frac{\partial R}{\partial\Theta}.
    \label{eq:pertubationterm}
\end{align}
This term is referred to as the pertubation term containing the design sensitivities in the following course of the work. Insertion of the derived expressions into Eq.~(\ref{eq:eulerlagrangedissipative}) yields subsequently the equation of motion
\begin{align*}
    \frac{\partial^2\phi}{\partial v^2} - c^2\nabla^2\phi + d\frac{\partial\phi}{\partial v} - h^4\nabla^4\phi = -f(\mathbf{u},\Theta;\mathbf{v},\lambda)\delta(\phi).
\end{align*}
In our research, we approximate the Heaviside function in Eq.~(\ref{eq:heavisidefunction}) using a hyperbolic tangent function
\begin{align}
    \Theta(\phi;\beta)\approx\frac{1}{2}(\tanh(2\beta\phi) + 1)
    \label{eq:hyperbolictangent}
\end{align}
where the parameter $\beta$ controls the width of the transition phase. Similar approaches were reported in \citep{PITSCH2008,COUPEZ2015,YAJI2016}. Differentiation of Eq.~(\ref{eq:hyperbolictangent}) with respect to $\phi$ yields the approximated Delta Dirac function as follows:
\begin{align}
    \delta(\phi;\beta)\approx\beta\mathrm{sech}^2(2\beta\phi).
    \label{eq:deltadiracapproximation}
\end{align}
Since we are only interested in the region of the transition phase between void and material domain, it is reasonable to reduce the complexity of Eq.~(\ref{eq:deltadiracapproximation}) by means of a Maclaurin series expansion, i.e.
\begin{align*}
    \beta\mathrm{sech}^2(2\beta\phi) = \sum_{n=0}^{\infty}\frac{(\beta\mathrm{sech}^2(2\beta\phi))^{(n)}(0)}{n!}\phi^n.
\end{align*}
Neglecting terms of higher order yields the simplified approximation of the Delta Dirac function
\begin{align}
    \delta(\phi;\beta) \approx \beta
    \label{eq:simplifieddeltadirac}
\end{align}
which then reduces to the simple constant $\beta$. By using the D'Alembert operator $\square(\cdot)\vcentcolon=(\partial^2/(c^2\partial v^2)-\nabla^2)(\cdot)$, considering Dirichlet and Neumann boundary conditions to be satisfied by $\phi$, and replacing the original expression of the Delta Dirac function by Eq.~(\ref{eq:simplifieddeltadirac}), we finally obtain the fundamental evolution problem
\begin{align}
    \left\{\begin{array}{rll}
    \displaystyle \square\phi + \frac{d}{c^2}\frac{\partial\phi}{\partial v} - \frac{h^4}{c^2}\nabla^4\phi &= \displaystyle -\frac{f(\mathbf{u},\Theta;\mathbf{v},\lambda)}{c^2}\beta  &\quad\text{in}\quad\mathscr{D}\\[10 pt]
    \phi &= 1 &\quad\text{on}\quad\mathscr{G}_\mathrm{t} \\[10 pt]
    \nabla\phi\cdot\mathbf{n} &=0 &\quad\text{on}\quad \partial\mathscr{D}\setminus\mathscr{G}_\mathrm{t} \\[10 pt]
    \phi\vert_{v=0} &= \phi_0 \\[10 pt]
    \displaystyle\left.\frac{\partial\phi}{\partial v}\right\vert_{v=0} &= \displaystyle\left(\frac{\partial\phi}{\partial v}\right)_0
    \end{array}\right.
    \label{eq:evolutionproblem}
\end{align} 
with initial conditions $\phi_0$ and $(\partial\phi/\partial v)_0$. In our work, the governing equation of motion is referred to as a damped generalized wave equation, as it incorporates a second derivative with respect to $v$ and both the Laplacian operator $\nabla^2(\cdot)$ and the biharmonic operator $\nabla^4(\cdot)$. As shown in Tab.~(\ref{tab:physicalequations}), our methodology offers significant flexibility, allowing us to model a wide range of fundamental physical phenomena by selectively considering different energetic contributions. Furthermore, the reaction-diffusion equation can also be derived as a special case from this framework, demonstrating that our approach aligns closely with the current state of research.
\begin{table}[ht]
    \centering
    \caption{Derivation of physical equations by combining different terms of energetic contribution.}
    \begin{adjustbox}{max width=\textwidth}
    \begin{tabular}{lccccc}\toprule
         \textsf{Type of equation} & $T[\partial\phi/\partial v]$ & $U_1[\phi]$ & $U_2[\phi]$ & $\mathscr{R}_1[\phi]$ & $\mathscr{R}_2[\partial\phi/\partial v]$ \\\midrule
         \textsf{Wave equation (WE)} & $\bullet$ & $\bullet$ & &$\bullet$& \\
         \textsf{Damped wave equation (DWE)} &$\bullet$ & $\bullet$ & &$\bullet$& $\bullet$  \\
         \textsf{Biharmonic wave equation (BWE)} & $\bullet$ & & $\bullet$ &$\bullet$&  \\
         \textsf{Damped biharmonic wave equation (DBWE)} & $\bullet$ & & $\bullet$ &$\bullet$& $\bullet$ \\
         \textsf{Generalized wave equation (GWE)} & $\bullet$ & $\bullet$ & $\bullet$ &$\bullet$&   \\
         \textsf{Damped Generalized wave equation (DGWE)} & $\bullet$ & $\bullet$ & $\bullet$ &$\bullet$& $\bullet$ \\\midrule
         \textsf{Reaction diffusion equation (RDE)} & & $\bullet$ &  &$\bullet$& $\bullet$  \\\bottomrule
    \end{tabular}
    \end{adjustbox}
    \label{tab:physicalequations}
\end{table}
Note that in case of the reaction-diffusion equation, the pertubation term $f(\mathbf{u},\Theta;\mathbf{v},\lambda)$ assumes the role of the reaction term, incorporating the design sensitivities. The according evolution problems given by equations~(\ref{eq:standardwaveequation}) to (\ref{eq:rdeequation}) and their weak formulations are collected in \ref{appendix}.

  \section{Numerical implementation}
\label{sec:numericalimplementation}
\begin{figure}
    \centering
    \includegraphics[width=1.\textwidth]{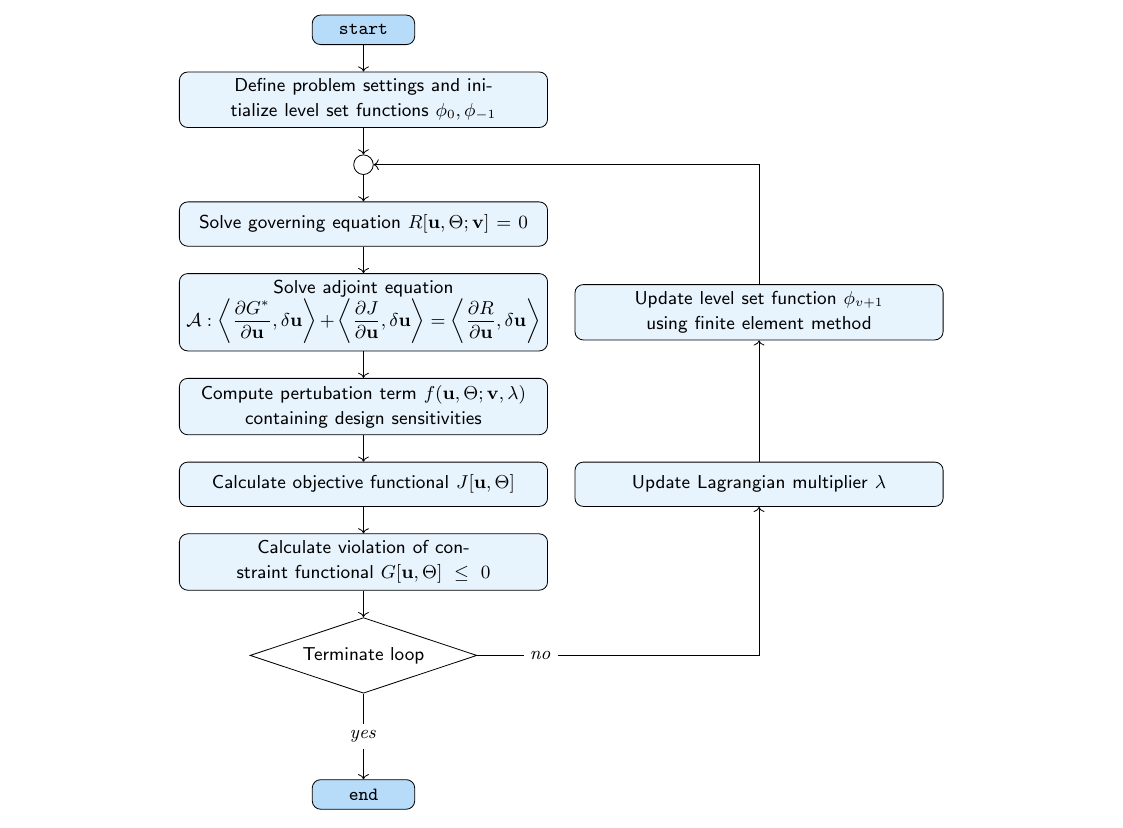}
    \caption{Flowchart to the optimization algorithm.}
    \label{fig:flowchart}
\end{figure}

We consider a discrete evolution step $\Delta v$ and apply the centered finite difference method to discretize the partial differential equation given in Eq.~(\ref{eq:evolutionproblem}) as follows:
\begin{align}
    \frac{\phi_{v+1} - 2\phi_{v} + \phi_{v-1}}{c^2\Delta v^2} + \frac{d}{c^2}\frac{\phi_{v+1} - \phi_v}{\Delta v} - \nabla^2\phi_{v+1} - \frac{h^4}{c^2}\nabla^4\phi_{v+1} = - \frac{f_v}{C_v\Delta v^2}\frac{\beta}{c^2}.
    \label{eq:finitedifferencemethod}
\end{align}
Here, the pertubation term $f_v$ is additionally normalized by a factor $C_v$ multiplied by $\Delta v^2$. In this context, we introduce the auxiliary parameters $\ell,m,k$ defined as $\ell = c\Delta v$, $m = d\Delta v$, and $k = h\sqrt{\Delta v}$ to further simplify Eq.~(\ref{eq:finitedifferencemethod}). After rearranging we obtain
\begin{align}
    (1 + m)\phi_{v+1} - \ell^2\nabla^2\phi_{v+1} - k^4\nabla^4\phi_{v+1} = - S_v + (2 + m)\phi_v - \phi_{v-1}
    \label{eq:discretizedequation}
\end{align}
while the normalized and discretized pertubation term is conveniently summarized to $S_v = f_v\beta / C_v$. The introduction of the auxiliary parameters has the advantage that the discrete evolutionary step $\Delta v$ no longer appears directly in the discredited equation but is indirectly taken into. This allows the individual effects of the parameters to be examined in isolation. Furthermore, parameters $\ell$ and $k$ then have the unit of a length $[\mathrm{L}]$ while $m$ remains dimensionless. In this connection, the presence of the biharmonic operator $\nabla^4(\cdot)$ in Eq.~(\ref{eq:discretizedequation}) poses a challenge in terms of solving the equation. Therefore, we introduce the substitution
\begin{align}
    \nabla^2\phi_{v+1} = \omega
    \label{eq:substitution}
\end{align}
and transform Eq.~(\ref{eq:discretizedequation}) into
\begin{align}
    (1+m)\phi_{v+1} - \ell^2\omega - k^4\nabla^2\omega = -S_v + (2+m)\phi_v - \phi_v.
    \label{eq:reducedequation}
\end{align}
This alternative formulation effectively reduces the highest derivative order, thereby simplifying the equation from a numerical perspective by lowering the continuity requirements for the basis functions in the subsequent finite element formulation. We now proceed to derive the weak form of Eq.~(\ref{eq:reducedequation}). Therefore, we first define the function space
\begin{align*}
    \Psi = \{\phi\in H^\alpha(\mathscr{D})\mid\phi=1\quad\text{on}\quad\mathscr{G}_\mathrm{t}\}
\end{align*}
and multiply Eq.~(\ref{eq:reducedequation}) by a suitable test function $\delta\phi\in\Psi$ followed by integration over $\mathscr{D}$:
\begin{align*}
    \int_\mathscr{D}(1+m)\phi_{v+1}\delta\phi\,\mathrm{d}\mathscr{D} - \int_\mathscr{D}\ell^2\omega\delta\phi\,\mathrm{d}\mathscr{D} - \int_\mathscr{D}k^4\nabla^2\omega\delta\phi\,\mathrm{d}\mathscr{D} = \dots \\ 
    \dots = \int_\mathscr{D}(-S_v + (2+m)\phi_v - \phi_{v-1})\delta\phi\,\mathrm{d}\mathscr{D}.
\end{align*}
For the substitution in Eq.~(\ref{eq:substitution}), we again multiply by a suitable test function $\delta\omega\in\Psi$ and obtain
\begin{align*}
    \int_\mathscr{D}\omega\delta\omega\,\mathrm{d}\mathscr{D} - \int_\mathscr{D}\nabla^2\phi_{v+1}\delta\omega\,\mathrm{d}\mathscr{D} = 0
\end{align*}
In the next step, we take the boundary conditions into account and apply the divergence theorem to transform the domain integrals in both equations. This leads to the coupled system of equations
\footnotesize{
\begin{align*}
    \left\{
    \begin{array}{rll}
        \displaystyle\int_\mathscr{D}(1+m)\phi_{v+1}\delta\phi\,\mathrm{d}\mathscr{D} - \int_\mathscr{D}\ell^2\omega\delta\phi\,\mathrm{d}\mathscr{D} + \int_\mathscr{D}k^4\nabla\omega\cdot\nabla\delta\phi\,\mathrm{d}\mathscr{D} &= \dots &\\[10pt]
        \displaystyle\dots = \int_\mathscr{D}(-S_v + (2+m)\phi_v - \phi_{v-1})\delta\phi\,\mathrm{d}\mathscr{D} & &\quad\text{in}\quad\mathscr{D}\\[10pt]
        \displaystyle\int_\mathscr{D}\omega\delta\omega\,\mathrm{d}\mathscr{D} + \int_\mathscr{D}\nabla\phi_{v+1}\cdot\nabla\delta\omega\,\mathrm{d}\mathscr{D} &= 0 &\quad\text{in}\quad\mathscr{D}\\[10 pt]
        \phi_{v+1} &= 1&\quad\text{on}\quad\mathscr{G}_\mathrm{t} \\[10 pt]
        \phi\vert_{v=0} &= \phi_0 & \\[10 pt]
        \displaystyle\left.\frac{\partial\phi}{\partial v}\right\vert_{v=0} &= \displaystyle\frac{\phi_0 - \phi_{-1}}{\Delta v} &
    \end{array}\right.
\end{align*}
}\normalsize
which can be solved using the finite element method. Therefore, the domain $\mathscr{D}$ is approximated by a set $\mathscr{E}$ consisting of finite elements $e$, so that
\begin{align*}
    \mathscr{D}\approx\bigcup_{e\in\mathscr{E}}e,
\end{align*}
while the functions are replaced with discretized vectors and matrices
\begin{align*}
    \mathbf{M} = \bigcup_{e\in\mathscr{E}}\int_{V_e}\mathbf{N}^\top\mathbf{N}\,\mathrm{d}V_e, && \mathbf{B} = \bigcup_{e\in\mathscr{E}}\int_{V_e}(\nabla\mathbf{N})^\top\nabla\mathbf{N}\,\mathrm{d}V_e
\end{align*}
containing the shape functions. Furthermore, we assume a generic evolution step of $\Delta v= 1$, therefore only the initial condition $\phi_{-1}$ has to be specified. This leads to the system of equations:
\footnotesize{
\begin{align*}
    \left\{
    \begin{array}{rll}
        (1 + m)\mathbf{M}\boldsymbol{\phi}_{v+1} - \ell^2\mathbf{M}\boldsymbol{\omega} + k^4\mathbf{B}\boldsymbol{\omega} &= \mathbf{M}(-\mathbf{s}_v + (2+m)\boldsymbol{\phi}_v - \boldsymbol{\phi}_{v-1}) &\quad\text{in}\quad\mathscr{E} \\[10 pt]
        \mathbf{M}\boldsymbol{\omega} + \mathbf{B}\boldsymbol{\phi}_{v+1} &= \mathbf{0} &\quad\text{in}\quad\mathscr{E} \\[10 pt]
        \boldsymbol{\phi}_{v+1} &= \mathbf{1} &\quad\text{on}\quad\bar{\mathscr{G}}_\mathrm{t} \\[10 pt]
        \boldsymbol{\phi}\vert_{v=0} &= \boldsymbol{\phi}_0 & \\[10 pt] 
        \boldsymbol{\phi}\vert_{v=-1} &= \boldsymbol{\phi}_{-1} &
    \end{array}\right.
\end{align*}
}\normalsize
Here, $\bar{\mathscr{G}}_\mathrm{t}$ denotes the discretized boundary segment $\mathscr{G}_\mathrm{t}
$. From the second equation we obtain the expression for $\boldsymbol{\omega}$ as follows:
\begin{align*}
    \boldsymbol{\omega} = -\mathbf{M}^{-1}\mathbf{B}\boldsymbol{\phi}_{v+1}.
\end{align*}
Insertion reduces the system of equations to a single equation given as
\footnotesize{
\begin{align}
    \left\{
    \begin{array}{rll}
        \left((1 + m)\mathbf{M} + \ell^2\mathbf{B} - k^4\mathbf{B}\mathbf{M}^{-1}\mathbf{B}\right)\boldsymbol{\phi}_{v+1} &= \mathbf{M}\left(-\mathbf{s}_v + (2+m)\boldsymbol{\phi}_{v} - \boldsymbol{\phi}_{v-1}\right) &\quad\text{in}\quad\mathscr{E} \\[10 pt]
        \boldsymbol{\phi}_{v+1} &= \mathbf{1} &\quad\text{on}\quad\bar{\mathscr{G}}_\mathrm{t} \\[10 pt]
        \boldsymbol{\phi}\vert_{v=0} &= \boldsymbol{\phi}_0 & \\[10 pt] 
        \boldsymbol{\phi}\vert_{v=-1} &= \boldsymbol{\phi}_{-1} &
    \end{array}\right.
\end{align}
}\normalsize
Figure~(\ref{fig:flowchart}) shows the corresponding flowchart of the optimization procedure. 
In the first step, level set functions $\phi_0$ and $\phi_{-1}$ are initialized, and we choose an order of $\alpha = 2$. Furthermore, additional parameters are specified, such as material properties. The governing equation is then solved using the finite element method. In the third step, the adjoint equation is computed to determine the pertubation term containing the design sensitivities. Based on this, the value of the objective functional and the constraint violation are evaluated. If convergence to a minimum is achieved, or a predefined number of iterations is reached, the loop is terminated. Otherwise, the Lagrange multiplier $\lambda$ and the level set function $\phi$ are updated, and the procedure is repeated. In this context, we use the software FreeFEM++ \citep{HECHT2012} to numerically solve the partial differential equations in combination with the mmg module \cite{DAPOGNY2014,BALARAC2022} which provides adaptive mesh refinements during the iterations.
  \section{Numerical experiments}
\label{sec:numericalexperiments}

\subsection{Minimum mean compliance}
Consider the benchmark problem of minimum mean compliance where a material domain $\Omega$ is fixed at $\Gamma_\mathrm{u}$ and loaded by a traction vector $\mathbf{t}$ at $\Gamma_\mathrm{t}$. For the sake of simplicity, body forces shall be neglected. Then, the according optimization problem reads 
\begin{customopti}|s|
    {inf}{\phi\in H^2(\mathscr{D})}{J[\mathbf{u}] = \int_{\Gamma_\mathrm{t}}\mathbf{t}\cdot\mathbf{u}\,\mathrm{d}\Gamma}{}{}{}
    \addConstraint{G[\Theta]=\frac{1}{V_0}\int_\mathrm{D}\Theta\,\mathrm{d}\Omega - V_\mathrm{f}}{\leq 0}{\quad\text{in}\quad\mathrm{D}}{}
    \addConstraint{-\mathrm{div}(\mathbf{C}:\varepsilon(\mathbf{u}))}{=0}{\quad\text{in}\quad\Omega}{}
    \addConstraint{\mathbf{u}}{=\mathbf{0}}{\quad\text{on}\quad\Gamma_\mathrm{u}}{}
    \addConstraint{(\mathbf{C}:\varepsilon(\mathbf{u}))\mathbf{n}}{=\mathbf{t}}{\quad\text{on}\quad\Gamma_\mathrm{t}}{}
    \addConstraint{(\mathbf{C}:\varepsilon(\mathbf{u}))\mathbf{n}}{=\mathbf{0}}{\quad\text{on}\quad\partial\Omega\backslash(\Gamma_\mathrm{u}\cup\Gamma_\mathrm{t})}{}
    \label{eq:minimummeancompliance}
\end{customopti}
with $V_\mathrm{f}$ being the desired volume fraction and $V_0$ the volume of the design domain, i.e.
\begin{align*}
    V_0 = \int_\mathrm{D}\,\mathrm{d}\Omega.
\end{align*}
From the boundary conditions, we derive the original governing equation with respect to $\Omega$ as
\begin{align*}
    \int_\Omega\varepsilon(\mathbf{u}):\mathbf{C}:\varepsilon(\mathbf{v})\,\mathrm{d}\Omega - \int_{\Gamma_\mathrm{t}}\mathbf{t}\cdot\mathbf{v}\,\mathrm{d}\Gamma = 0\quad\mathbf{u},\mathbf{v}\in\mathscr{U}
\end{align*}
in its weak form. To extend the integral equation to the entire design domain $\mathrm{D}$, we relate to an Ersatz material approach \citep{CHOI2011,WANG2003,ALLAIRE2004} to express the elasticity tensor depending on $\phi$ by introducing the auxiliary function
\begin{align}
    \tau(\Theta) = (1-e)\Theta^q(\phi) + e 
    \label{eq:tau}
\end{align}
with $e\ll 1$ and $q>1$ so that it takes the form
\begin{align*}
    \mathbf{C}_\phi(\Theta) = \tau(\Theta)\mathbf{C}\quad\text{in}\quad\mathrm{D}.
\end{align*}
Consequently, the partial derivative of Eq.~(\ref{eq:tau}) with respect to $\Theta$ results to
\begin{align*}
    \frac{\partial\tau}{\partial\Theta} = q(1-e)\Theta^{q-1}(\phi).
\end{align*}
Since $\Theta(\phi)\in\{0,1\}$ per definition in Eq.~(\ref{eq:heavisidefunction}), raising to any power $q>1$ does not change its value. Therefore, we can use the identities
\begin{align*}
    \Theta^q(\phi) \equiv \Theta^{q-1}(\phi) \equiv \Theta(\phi)
\end{align*}
to simplify the expression so that
\begin{align}
    \frac{\partial\tau}{\partial\Theta} = q(1 - e)\Theta(\phi)
    \label{eq:dtau}
\end{align}
holds. Then, we obtain
\begin{align*}
    R_1:\quad\int_\mathrm{D}\varepsilon(\mathbf{u}):\tau(\Theta)\mathbf{C}:\varepsilon(\mathbf{v})\,\mathrm{d}\Omega = \int_{\Gamma_\mathrm{t}}\mathbf{t}\cdot\mathbf{v}\,\mathrm{d}\Gamma 
\end{align*}
for the governing equation in level set formulation. The partial functional derivative of the constraint functional with respect to the displacement field vanishes, therefore the adjoint equation is governed from the remaining terms and given as 
\begin{align*}
    \mathcal{A}_1:\quad\int_{\Gamma_\mathrm{t}}\mathbf{t}\cdot\delta\mathbf{u}\,\mathrm{d}\Gamma = \int_\mathrm{D}\varepsilon(\delta\mathbf{u}):\tau(\Theta)\mathbf{C}:\varepsilon(\mathbf{v})\,\mathrm{d}\Omega
\end{align*}
according to Eq.~(\ref{eq:adjointequation}). Since the considered problem is self-adjoint, it directly follows from the adjoint equation that the state and adjoint variables coincide, i.e., $\mathbf{u} = \mathbf{v}$. In our research, the constraint functional is expressed by means of the Augmented Lagrangian method so that
\begin{align}
   G^\ast[\Theta] = \lambda\bar{G}[\Theta] + \frac{r}{2}\bar{G}^2[\Theta] 
\end{align}
applies while the constraint functional is defined as
\begin{align*}
    \bar{G}[\Theta] = \max\left\{G[\Theta];-\frac{\lambda^\prime}{r}\right\}.
\end{align*}
In this context, $\lambda^\prime$ denotes the multiplier of the previous iteration step and $r$ represents an amplification factor which is increased in each iteration until it reaches a limit value $r_\mathrm{max}$. Thus, the partial functional derivative of the constraint functional with respect to $\Theta$ follows to 
\begin{align*}
    \left\langle\frac{\partial G^\ast}{\partial\Theta},\delta\Theta\right\rangle &= \frac{\mathrm{d}}{\mathrm{d}\epsilon}\left[\lambda\bar{G}[\Theta + \epsilon\delta\Theta] + \frac{r}{2}\bar{G}^2[\Theta + \epsilon\delta\Theta]\right]_{\epsilon=0} \\
    &= (\lambda + r\bar{G}[\Theta])\int_\mathrm{D}\frac{\partial\bar{g}}{\partial\Theta}\delta\Theta\,\mathrm{d}\Omega
\end{align*}
while the partial derivative in the integral results to
\begin{align*}
    \frac{\partial\bar{g}}{\partial\Theta} = 
    \left\{\begin{array}{ll}
        \displaystyle \frac{1}{V_0} &\displaystyle\quad\text{if}\quad G[\Theta]>-\frac{\lambda^\prime}{r} \\[10 pt]
        0 &\quad\text{otherwise}
    \end{array}\right..
\end{align*}
Since this is valid for all test functions $\delta\Theta$, the derivative reads
\begin{align*}
    \frac{\partial G^\ast}{\partial\Theta} = (\lambda + r\bar{G}[\Theta])\frac{\partial\bar{g}}{\partial\Theta}.
\end{align*}
Following the update rule of the Augmented Lagrangian method, the multiplier is updated via
\begin{align*}
    \lambda = \lambda^\prime + r\max\{G[\Theta];0\}.
\end{align*}
In this context, we conveniently define an indicator function
\begin{align*}
    \mathbb{I}_G \vcentcolon =
    \left\{\begin{array}{ll}
        1 &\quad\text{if}\quad \lambda^\prime + r G[\Theta] > 0 \\ [10 pt]
        0 &\quad\text{otherwise} 
    \end{array}\right.,
\end{align*}
which allows the single augmented terms to be reformulated as follows:
\begin{align*}
    \bar{G}[\Theta] &= G[\Theta]\mathbb{I}_G + \left(-\frac{\lambda^\prime}{r}\right)(1 - \mathbb{I}_G), \\
    \frac{\partial\bar{g}}{\partial\Theta} &= \frac{1}{V_0}\mathbb{I}_G.
\end{align*}
Insertion then first yields
\begin{align}
    \frac{\partial G^\ast}{\partial\Theta} = (\lambda + r G[\Theta]\mathbb{I}_g - \lambda^\prime(1 - \mathbb{I}_g))\frac{1}{V_0}\mathbb{I}_G.
    \label{eq:partialderivativevolumeconstraint}
\end{align}
Here, we use the identity $\mathbb{I}_G^2\equiv\mathbb{I}_G$ and obtain
\begin{align*}
    \frac{\partial G^\ast}{\partial\Theta} = (\lambda + r G[\Theta])\frac{1}{V_0}\mathbb{I}_G
\end{align*}
for the partial derivative of the constraint functional. According to Eq.~(\ref{eq:pertubationterm}) and using Eq.~(\ref{eq:dtau}), we finally obtain the pertubation term as
\begin{align*}
    f_1 = (\lambda + r G[\Theta])\frac{1}{V_0}\mathbb{I}_G - \varepsilon(\mathbf{u}):\frac{\partial\tau}{\partial\Theta}\mathbf{C}:\varepsilon(\mathbf{u})
\end{align*}
and additionally formulate the normalization factor
\begin{align}
    C_v = \frac{1}{\mathrm{c}_fV_0}\int_\mathrm{D}\left\vert\left(\frac{\partial R}{\partial\Theta}\right)\right\vert_{L^2(\mathrm{D})}\,\mathrm{d}\Omega
    \label{eq:standardnormalizationfactor}
\end{align}
with constant $\mathrm{c}_f>0$. We have now established all equations required for numerically solving the evolution problem, namely the equilibrium condition, the adjoint equation, and the perturbation term.

\begin{figure}[ht]
    \centering
    \includegraphics[width=1.\textwidth]{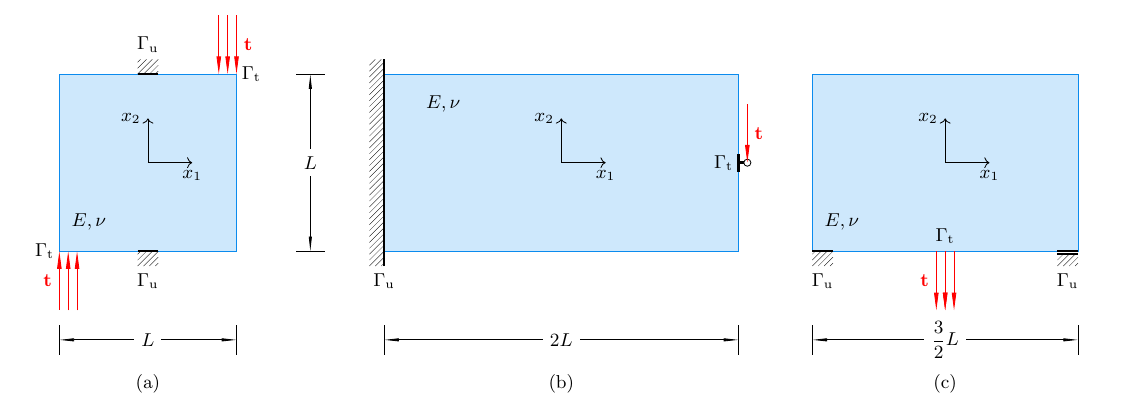}
    \caption{Different two-dimensional problems; (a) asymmetrically loaded plate, (b) cantilever beam, (c) simply supported girder.}
    \label{fig:testcases}
\end{figure}

In the following, the problem of minimum mean compliance shall serve to investigate the impact of the individual numerical parameters on the example of different two-dimensional problem settings given in Fig.~(\ref{fig:testcases}). 
In all scenarios, the admissible volume fraction is $V_\mathrm{f} = 0.45$ while Young's modulus of the isotropic linear elastic material is set to $E = 210\,\mathrm{GPa}$, Poisson's ratio to $\nu=0.3$ and the traction vector's magnitude to $|\mathbf{t}|=1.0\times 10^3\,\mathrm{N/m}$. The length is fixed at $L=1.0\,\mathrm{m}$. In addition, the parameters corresponding to the Ersatz material approach are set to $q = 3.0$ and $e = 1.0\times 10^{-3}$. The finite element size ranges from $h_\mathrm{min} = 5\times 10^{-4}$ to $h_\mathrm{max} = 1.0\times 10^{-2}$. The proportionality constant is set to $\mathrm{c}_f = 1.0$ and a completely filled design domain serves as the initial configuration for both $\phi_0$ and $\phi_{-1}$.

\subsubsection{Impact of the width of transition phase $\beta$}
We consider an asymmetrically loaded square plate being fixed at $\Gamma_\mathrm{u}$ as shown in Fig.~(\ref{fig:testcases}a) and use the update scheme of the standard wave equation with parameters $\ell = 0.008, m=0.0, k=0.0$. The width of the transition phase is varied between $\beta=2.0$ and $\beta=8.0$.  
\begin{figure}[ht]
    \centering
    \includegraphics[width=1.\textwidth]{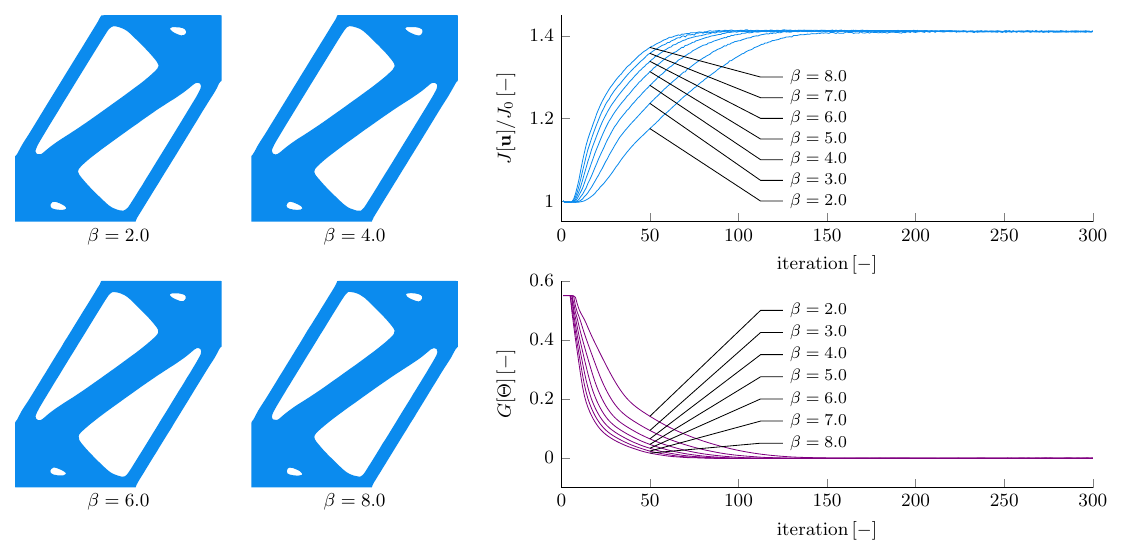}
    \caption{Influence of transition phase parameter $\beta$.}
    \label{fig:evaluationbeta}
\end{figure}
In this context, the resulting topologies are as well as the convergence histories of both the ratio of the objective functional $J[\mathbf{u}]/J_0$ and the constraint functional $G[\Theta]$ are depicted in Fig.~(\ref{fig:evaluationbeta}). Here, it can be seen that the topologies remain almost identical with varying $\beta$. Nevertheless, we further notice that an increase in $\beta$ increases the speed of convergence. This observation can be explained by the fact that $\beta$ proportionally amplifies the perturbation term which acts as a driving force. This further suggests that the external excitation of the fictitious matter intensifies proportionally with the amount of supplied energy.

\subsubsection{Impact of parameters on the wave equation}
\begin{figure}[ht]
    \centering
    \includegraphics[width=1.\textwidth]{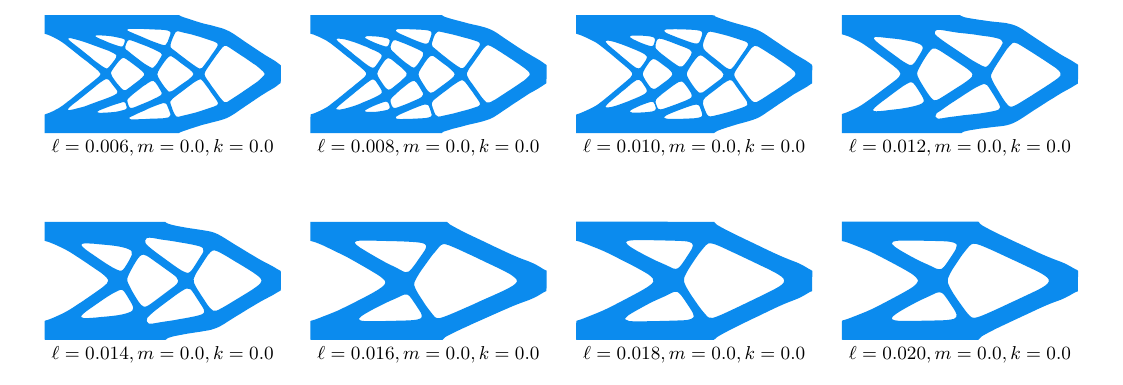}
    \caption{Influence of parameter $\ell$ on the topological layout in the standard wave equation.}
    \label{fig:evaluationwave}
\end{figure}

In the following section, the influences of the wave parameters will be examined in more detail using the cantilever beam problem from Fig.~(\ref{fig:testcases}b), whereby the transition phase parameter is set to $\beta = 5.0$ for all experiments. We start with the standard wave equation and vary the parameter $\ell$ between $0.006$ and $0.020$. Here, Fig.~(\ref{fig:evaluationwave}) shows the resulting topologies for different $\ell$.
We find that the parameter $\ell$ in the standard wave equation has a direct impact on the complexity of the structural layout. In this context, the complexity characterizes the number of holes created within the structure. It can be seen that small values of $\ell$ result in a significantly more complex topology, while an increase in $\ell$ leads to a more robust layout, i.e. the number of holes decreases. This effect can be attributed to the fact that $\ell$ governs the resistance of the fictitious matter to both compression and tension through parameter $c$, see Fig.~(\ref{fig:transition}). Specifically, at higher values of $\ell$ the material becomes more resistant to deformation, resulting in reduced strain and, consequently, fewer holes being formed.

We now wish to examine how additional damping affects the results and apply the update scheme of the damped wave equation for this purpose. In this context, it is useful to classify the resulting topologies into the categories "complex", "moderate" and "robust". This classification serves as a qualitative measure to describe the number of holes contained. We apply a similar measure to the damping and distinguish between "weak", "medium" and "strong" damping.
\begin{figure}[ht]
    \centering
    \includegraphics[width=1.\textwidth]{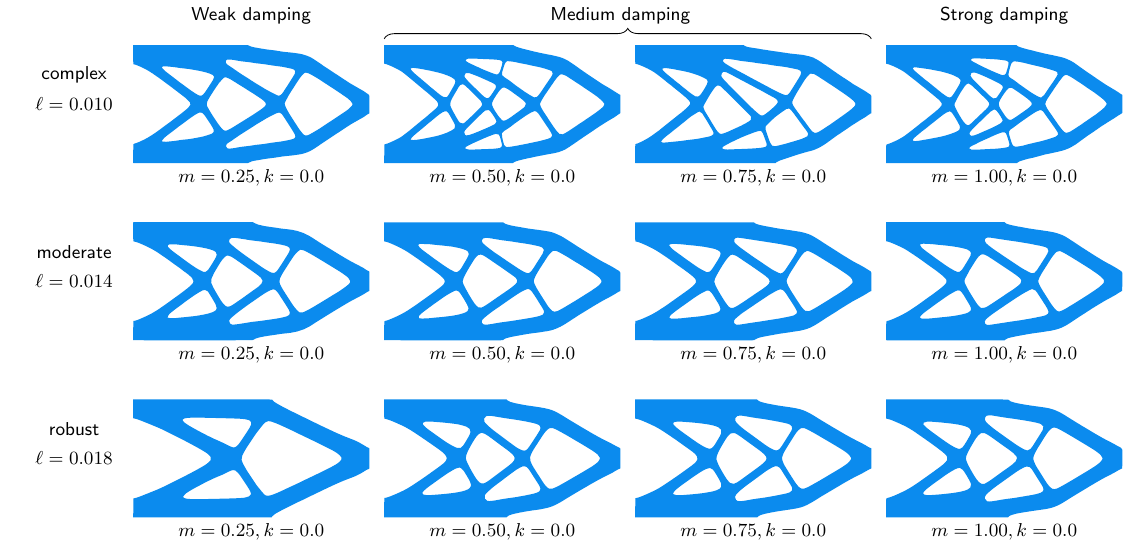}
    \caption{Influence of parameter $m$ on the structural layout in the damped wave equation.}
    \label{fig:evaluationLm}
\end{figure}
Figure~(\ref{fig:evaluationLm}) shows the resulting topologies for different levels of complexity and damping. It can be seen, that in the case of moderately complex topologies, damping has only a minor influence. In contrast, for the robust variant, increasing damping leads to higher complexity, gradually transitioning toward a moderately complex layout. For initially more complex topologies, weak damping results in a moderate layout, whereas increasing the damping again leads to higher complexity and heterogeneous patterns.
Notably, the moderate topology appears to represent a preferred solution, forming a pronounced plateau in the search space that is preferably assumed. Overall, the moderate layout emerges across all levels of complexity, suggesting that the parameter $m$ can serve as an indicator for identifying this favored configuration.
\begin{figure}[ht]
    \centering
    \includegraphics[width=1.\textwidth]{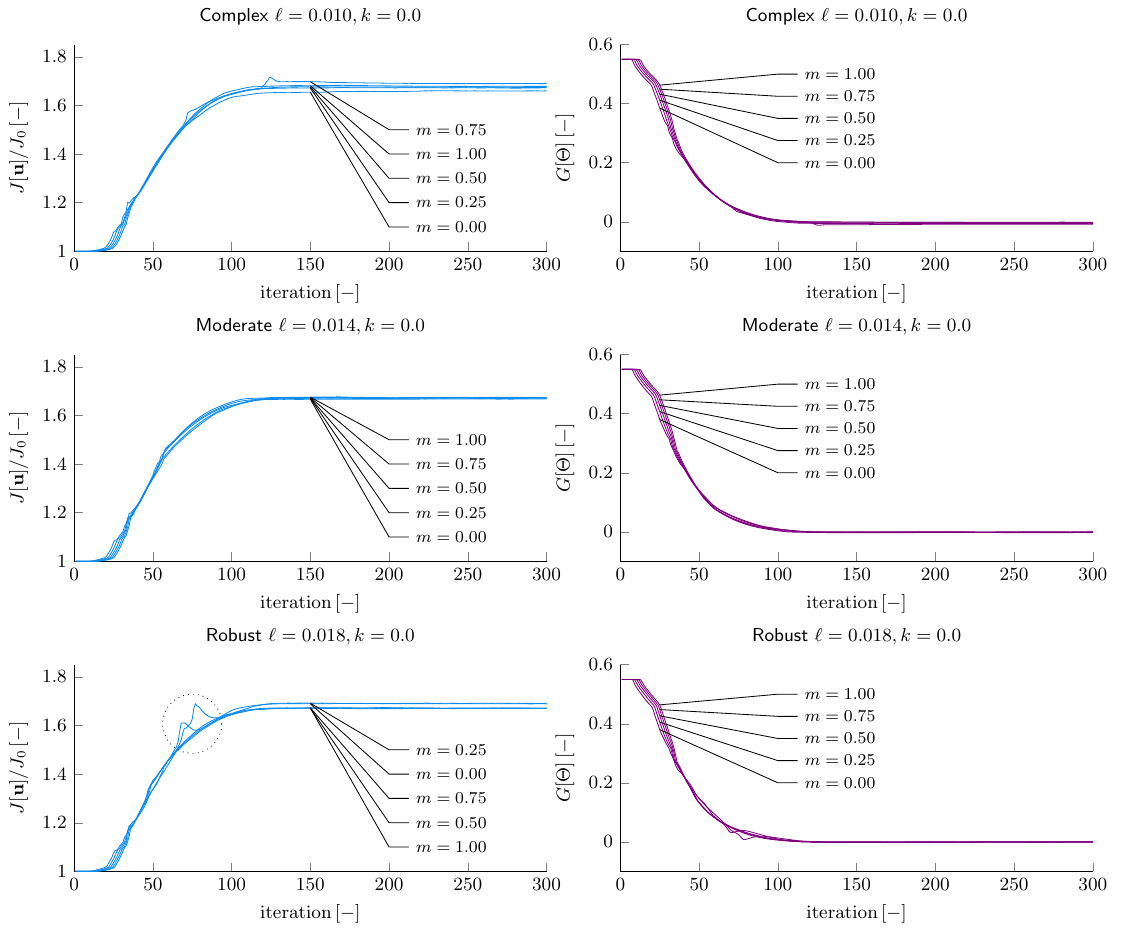}
    \caption{Influence of parameter $m$ on the convergence behavior in the damped wave equation.}
    \label{fig:evaluationLmdiagram}
\end{figure}
The convergence curves shown in Fig.~(\ref{fig:evaluationLmdiagram}) also indicate that the convergence speed is slightly reduced with increasing damping. The pronounced deflections (see for instance the encircled region in the diagram of the robust variant) mark areas where struts are removed, which influences the compliance in the short term before the structure stabilizes again in its new configuration.

\subsubsection{Impact of parameters on the biharmonic wave equation}
\begin{figure}[ht]
    \centering
    \includegraphics[width=1.\textwidth]{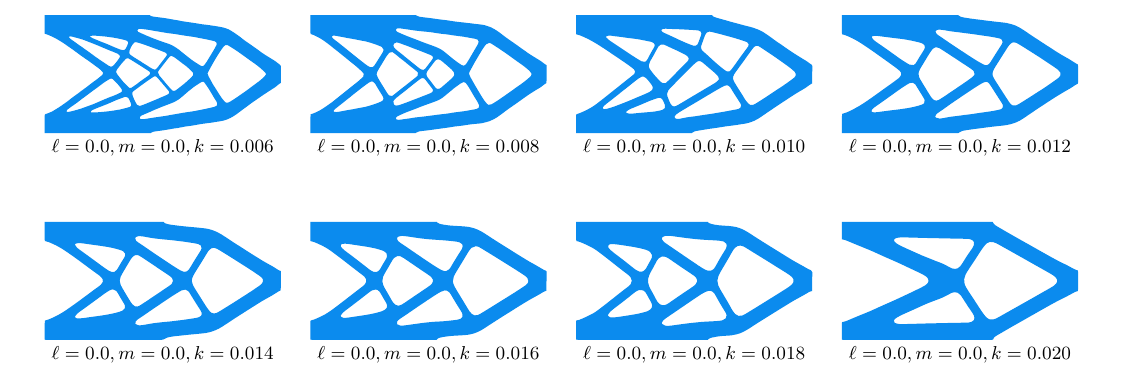}
    \caption{Influence of parameter $k$ on the structural layout in the biharmonic wave equation.}
    \label{fig:evaluationbiharmonic}
\end{figure}

\begin{figure}[ht]
    \centering
    \includegraphics[width=1.\textwidth]{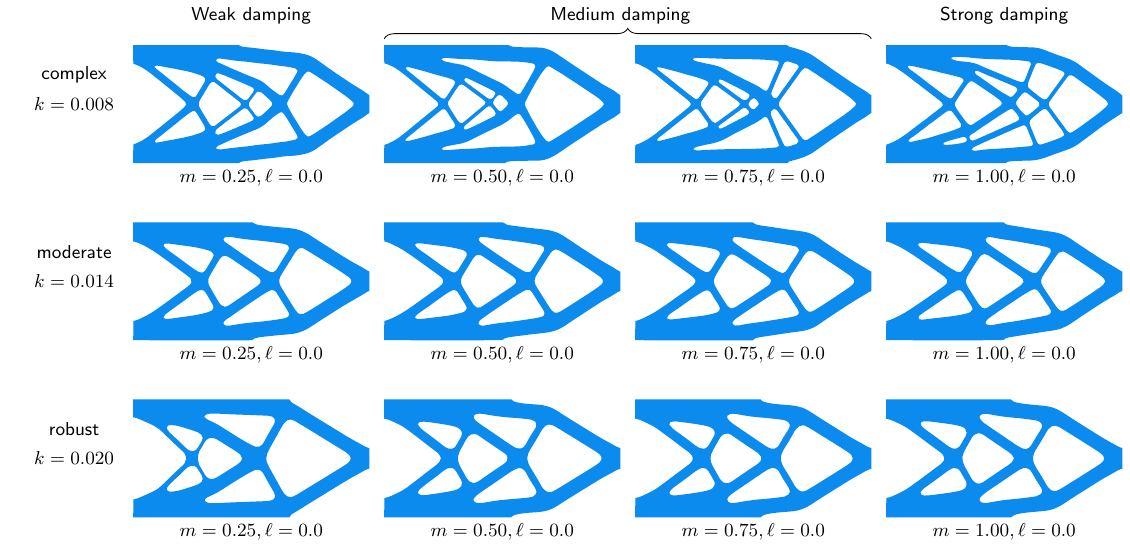}
    \caption{Influence of parameter $m$ on the structural layout in the damped biharmonic wave equation.}
    \label{fig:evaluationdampedbiharmonic}
\end{figure}

\begin{figure}[ht]
    \centering
    \includegraphics[width=1.\textwidth]{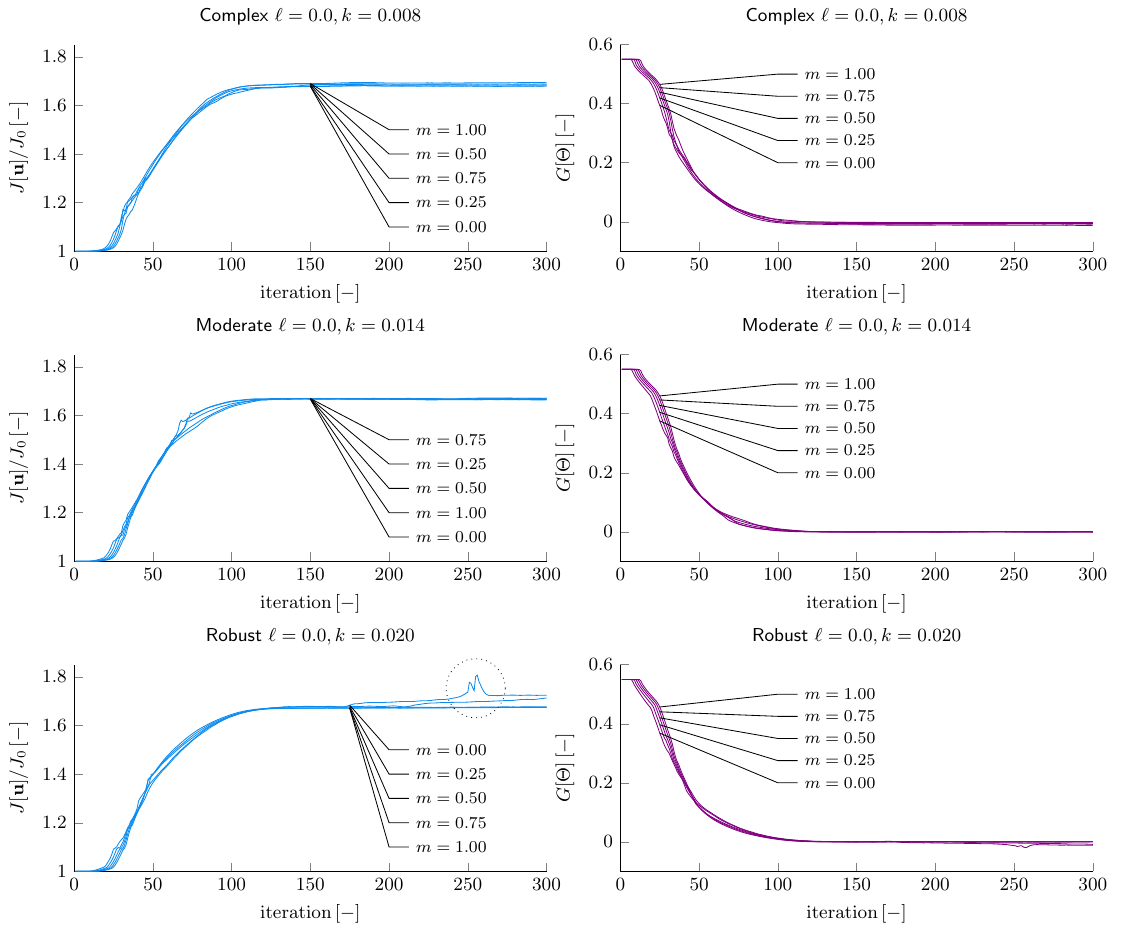}
    \caption{Influence of parameter $m$ on the convergence behavior in the damped biharmonic wave equation.}
    \label{fig:evaluationKmdiagram}
\end{figure}

With regard to the biharmonic wave equation, a similar pattern can be seen in the increase in topological complexity. As the parameter $k$ increases, the number of holes decreases, resulting in a more robust structure. While the development of strut thickness appears almost linear in the case of the standard wave equation for complex topologies, a more nonlinear pattern can be observed in the biharmonic case (cf. variants $k = 0.006$ and $k = 0.008$ in Fig.~(\ref{fig:evaluationbiharmonic})). For instance, fine interior structures are formed, which are enclosed by more robust outer struts. However, the influence of this effect decreases with increasing $k$.
Furthermore, we analyze the influence of additional damping on the resulting topologies of different levels of complexity. As illustrated in Fig.~(\ref{fig:evaluationdampedbiharmonic}), damping has only a minor impact on the number of holes, but significantly affects the geometric layout of the structures. This is particularly evident in the complex variant, where the overall degree of topological complexity remains almost constant, however, the structural layout changes significantly. For instance, at $m=0.75$, the front struts tend to split, and the created holes widen with increasing damping. In contrast, the moderate and robust variants exhibit behavior similar to that observed in the damped standard wave equation: the system consistently favors the moderate topology, which appears to emerge as a preferred solution.
With regard to the convergence histories shown in Fig.~(\ref{fig:evaluationKmdiagram}), we also observe with the damped biharmonic wave equation that the speed of convergence is reduced by the damping. The peaks (such as in the diagram of the robust variant) in the graph are also due to the removal of structural elements.

\subsubsection{Impact of parameters on the generalized wave equation}
\begin{figure}[ht]
    \centering
    \includegraphics[width=1.\textwidth]{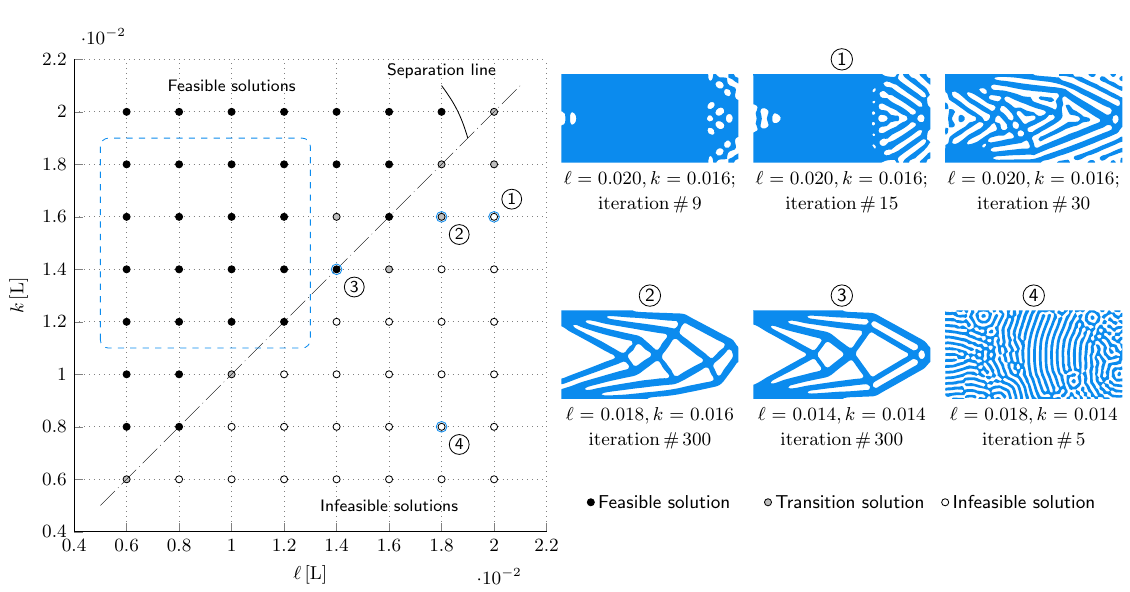}
    \caption{Influence of combined parameters $\ell$ and $k$ on the feasibility of solutions in the generalized wave equation.}
    \label{fig:evaluationlkdiagram}
\end{figure}

Finally, we analyze the influence of the parameters $k$ and $\ell$ in the generalized wave equation. The numerical experiments show that not every combination of these parameters leads to a useful structural solution. In particular, pronounced wave patterns can be observed with unfavorable parameter settings, which dominate and prevent structure formation (see configuration "4" in Fig.~(\ref{fig:evaluationlkdiagram}) and configuration series "1"). In case of configuration "1", the structure remains partially recognizable, but lacks clear functional suitability. Clearly defined, functional topologies only emerge from certain parameter combinations, as shown as an example in configuration "3" of Fig.~(\ref{fig:evaluationlkdiagram}). A detailed analysis of the parameter combinations suggests that the limit between functionally stable and non-stable solutions can be located along a characteristic separation line in the $\ell-k$ parameter space, see left diagram in Fig.~(\ref{fig:evaluationlkdiagram}). This dashed depicted trend line runs approximately along the ratio $\ell/k \approx 1$. It can therefore be derived from the diagram that the ratio $\ell/k<1$ should be maintained for the emergence of suitable topologies. Close to the separation line, the wave-like patterns increase significantly, causing the structures to lose coherence and begin to disintegrate as it can be seen in configuration "2" of Fig.~(\ref{fig:evaluationlkdiagram}).

\begin{figure}[ht]
    \centering
    \includegraphics[width=1.\textwidth]{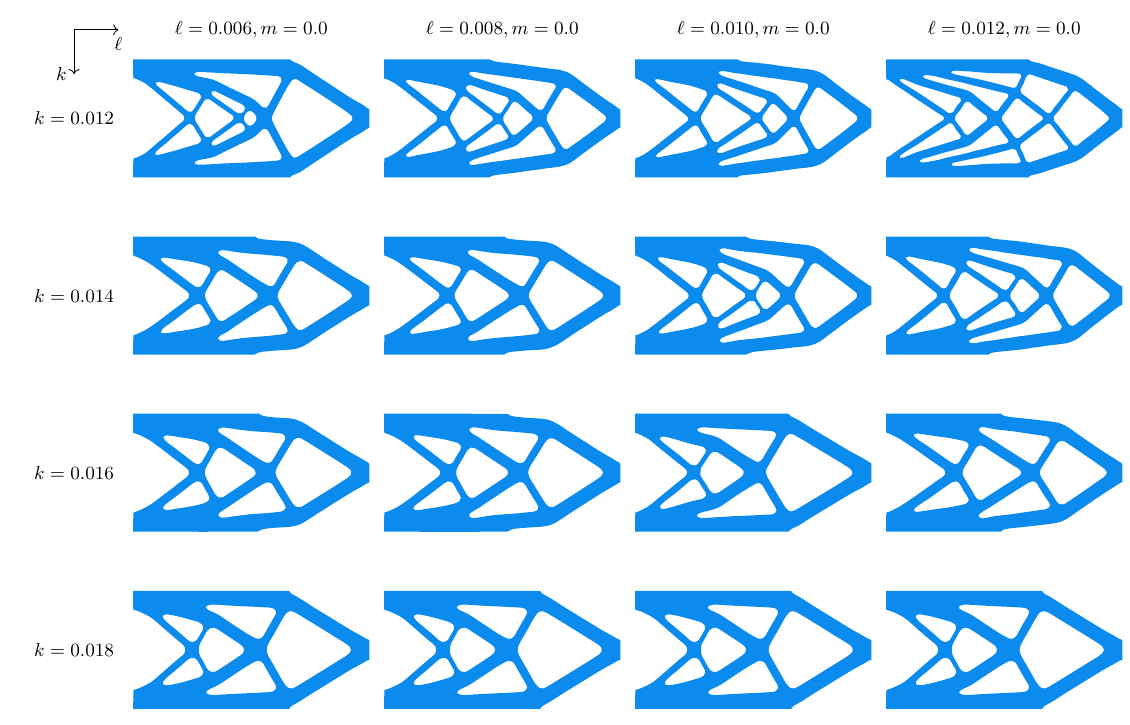}
    \caption{Influence of combined parameters $\ell$ and $k$ on the structural layout in the generalized wave equation.}    
    \label{fig:evaluationlk}
\end{figure}

\begin{figure}[ht]
    \centering
    \includegraphics[width=1.\textwidth]{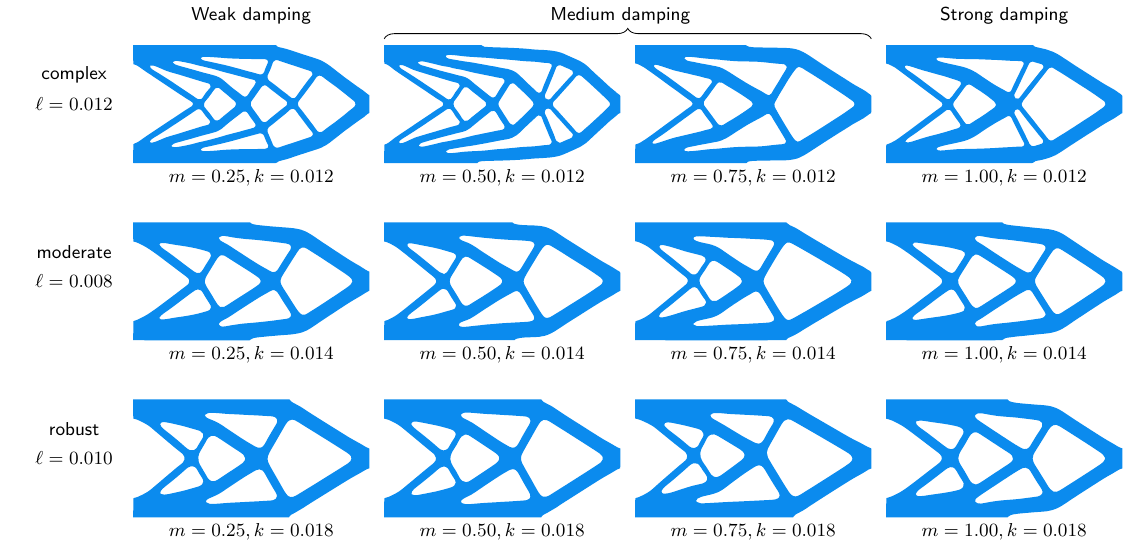}
    \caption{Influence of parameter $m$ on the structural layout in the damped generalized wave equation.}
    \label{fig:evaluationLmk}
\end{figure}

\begin{figure}[ht]
    \centering
    \includegraphics[width=1.\textwidth]{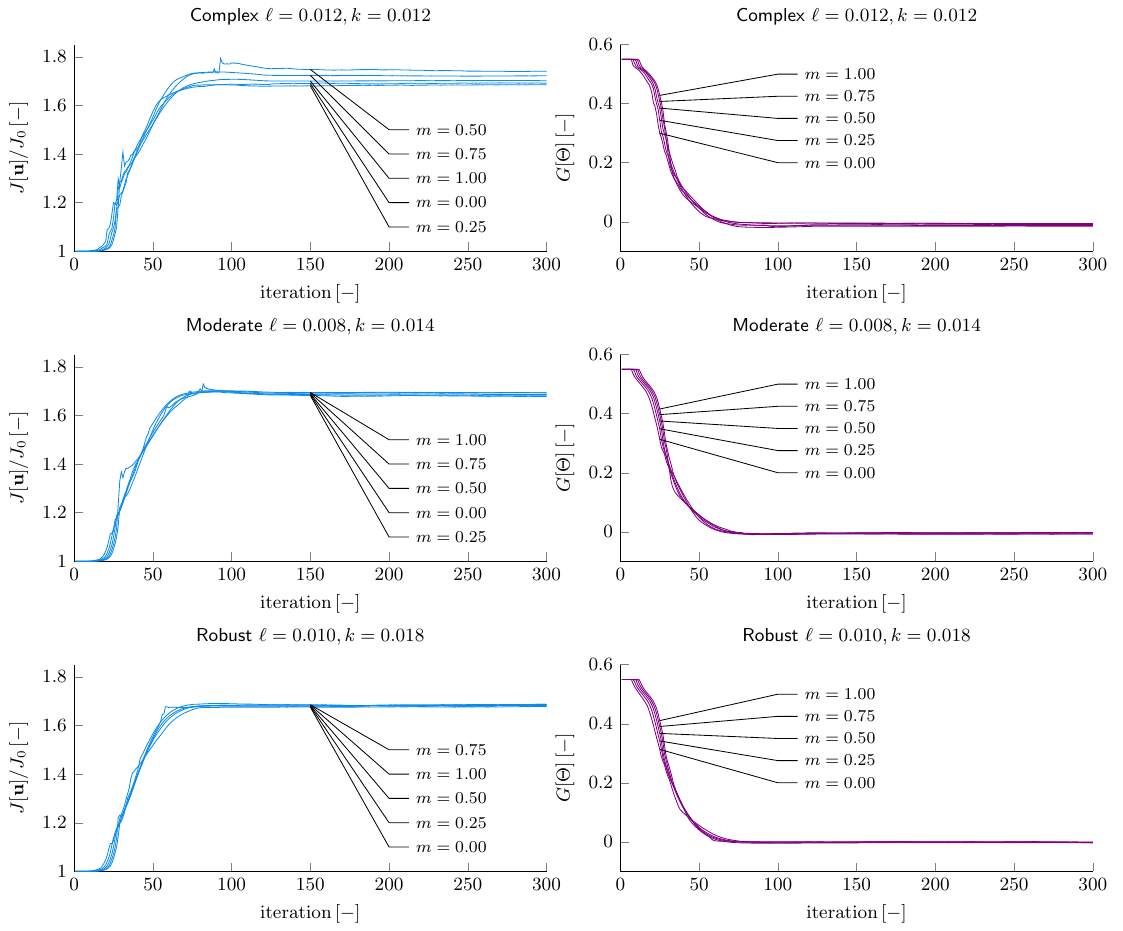}
    \caption{Influence of parameter $m$ on the convergence behavior in the damped generalized wave equation.}
    \label{fig:evaluationLmkdiagram}
\end{figure}

Furthermore, a remarkable effect can be observed with regard to the geometry of the struts: These have a bulbous shape in the middle, while tapering towards the connecting points, see for instance the shape of the inner struts of configuration "2" or "3" in Fig.~(\ref{fig:evaluationlkdiagram}). This effect becomes more pronounced near the separation line. Such a pattern is reminiscent of classic column shapes, where the column slightly curves towards the center. This principle can be found both in architecture and in biological structures and in both cases serves to distribute material efficiently under load. In the context of topology optimization, the emergent strut shape indicates an adaptation to local stress distributions.
Based on this, the resulting topologies within the section marked with a dashed line in Fig.~(\ref{fig:evaluationlkdiagram}) are analyzed in more detail in Fig.~(\ref{fig:evaluationlk}). It can be seen that for small values of $k$, an increase in the parameter $\ell$ influences the complexity of the structure: existing holes become larger and the struts tend to split. This effect decreases significantly with increasing $k$. In this range, it can be observed that a moderate topology is preferred.
Finally, the influence of additional damping in the generalized wave equation will also be investigated in more detail, see Fig.~(\ref{fig:evaluationLmk}). In the case of a complex topology, increasing damping leads to a slight reduction in complexity. In contrast, the moderate and robust variant shows only a slight influence of damping. However, it is noticeable that a moderate layout is preferred in both cases.
With regard to the convergence curves, a similar conclusion can be drawn as with the curves of the damped biharmonic and damped standard wave equation. Here too, damping reduces the speed of convergence as shown in Fig.~(\ref{fig:evaluationLmkdiagram}). 

\subsubsection{Symmetry aspects}
When analyzing the topologies obtained in the numerical experiments, it is noticeable that asymmetrical patterns arise under certain parameter settings - even though a symmetrical solution would be expected due to the symmetrical boundary conditions (as in the cantilever beam example). This behavior can be explained by the fact that the system remains in a locally stable but asymmetrical solution. In these cases, the driving force induced by the perturbation term is not sufficient to transfer the system into a symmetric configuration. A central aspect of this behavior is the direct coupling between the transition phase parameter $\beta$ and the perturbation term, which can be used specifically to control the intensity of the driving force. A reduction of the parameter $\beta$ reduces this force and can lead the system into a neighboring, but symmetric equilibrium position. An increase in $\beta$, on the other hand, increases the driving force, resulting in symmetrical and more robust topologies. 
\begin{figure}[ht]
    \centering
    \includegraphics[width=1.\textwidth]{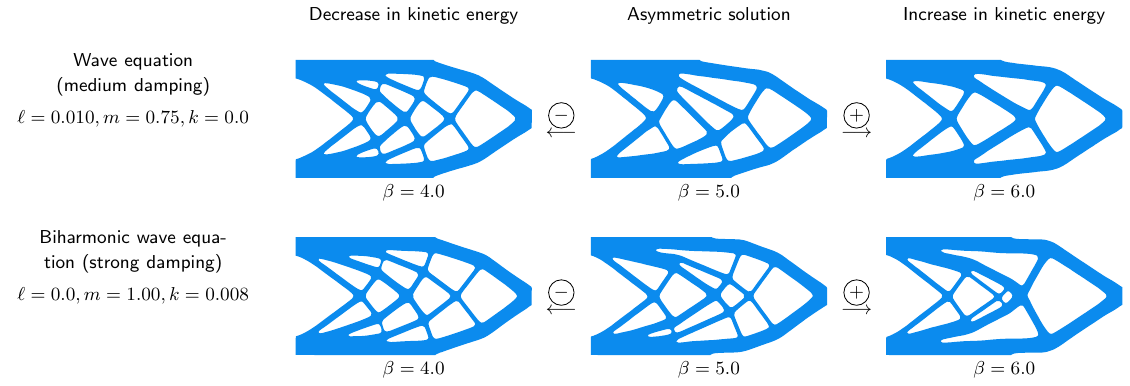}
    \caption{Transition to symmetrical layout by varying parameter $\beta$.}
    \label{fig:symmetry}
\end{figure}
We demonstrate this behavior using two asymmetric solutions shown in Fig.~(\ref{fig:symmetry}) in which $\beta$ is specifically varied: By increasing or reducing the parameter, symmetric solutions can be generated in both cases. It is also shown that the structural complexity increases with decreasing $\beta$, while it decreases with increasing $\beta$. This behavior is in line with our interpretation formulated earlier.

\subsubsection{Influence of initial configuration}
\begin{figure}[ht]
    \centering
    \includegraphics[width=1.\textwidth]{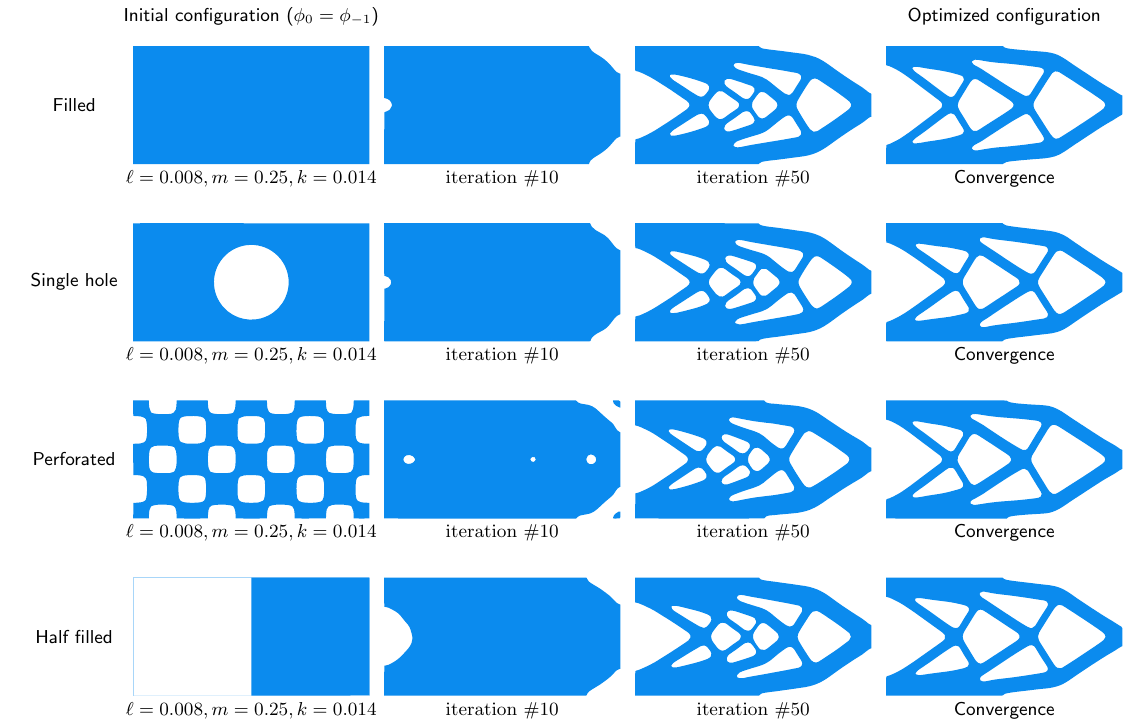}
    \caption{Influence of different initial configurations $\phi_0=\phi_{-1}$ on the solution considering the damped generalized wave equation.}
    \label{fig:initialdampedgeneralized}
\end{figure}

\begin{figure}[ht]
    \centering
    \includegraphics[width=1.\textwidth]{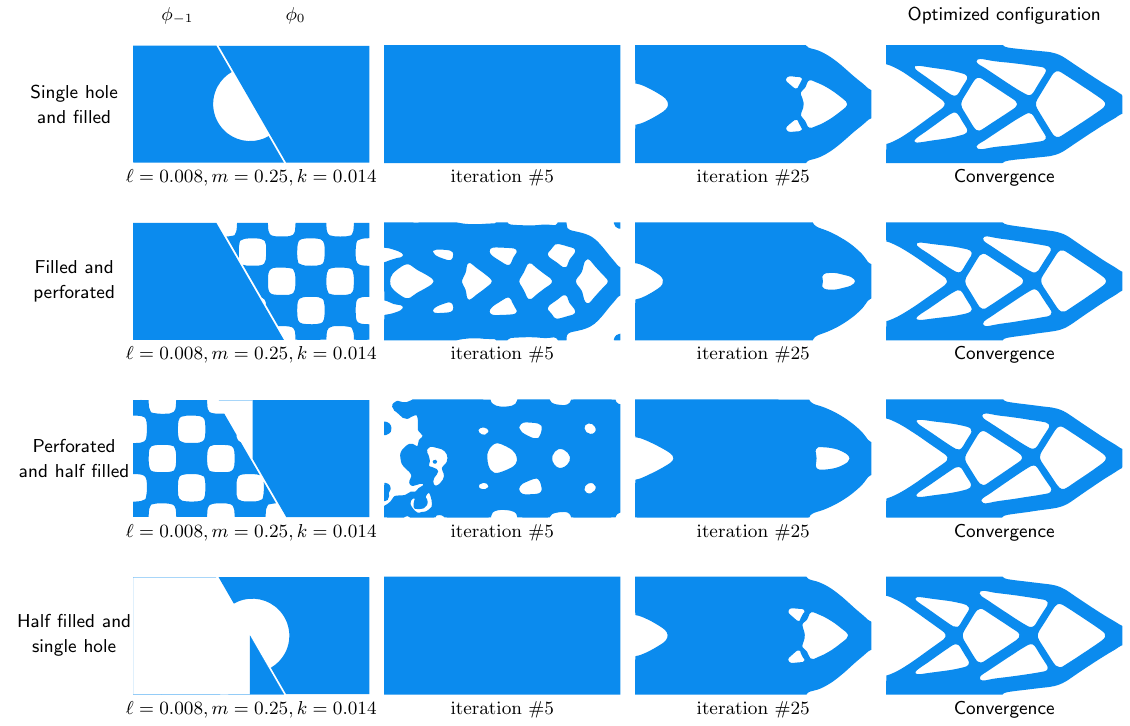}
    \caption{Influence of different combined initial configurations $\phi_0$ and $\phi_{-1}$ on the solution considering the damped generalized wave equation.}
    \label{fig:initialmixeddampedgeneralized}
\end{figure}

To assess the robustness of the proposed method, we investigate the influence of different initial configurations, defined by the level set functions $\phi_0$ and $\phi_{-1}$, on the final topology. As a reference, we employ the weakly damped generalized wave equation with $\beta=5.0$, which incorporates all model parameters. The remaining numerical parameters are kept unchanged throughout this study. In a first step, we assume $\phi_0=\phi_{-1}$ and use a completely filled solid design domain as the baseline configuration. We then analyze the effect of three perturbations: a single hole, a perforated material distribution, and a half-filled domain. The resulting topologies and their evolution across iterations are shown in Fig~(\ref{fig:initialdampedgeneralized}). 
In all cases, the same optimized configuration is obtained, demonstrating insensitivity to the initial configuration. Notably, the initial disturbances are rapidly dissipated: after approximately 10 iterations, the solution begins to stabilize, and after 50 iterations, the intermediate shapes are identical. Furthermore, we examine the influence of mixed initial configurations, i.e. $\phi_0\neq\phi_{-1}$, combining different levels of material and perforation, see Fig.~(\ref{fig:initialmixeddampedgeneralized}). 
Again, the algorithm converges to the same final structural layout, confirming that the method is robust with respect to both spatial and temporal initialization.

\subsubsection{Comparison of the different physical equations of motion}
\begin{figure}[ht]
    \centering
    \includegraphics[width=1.0\textwidth]{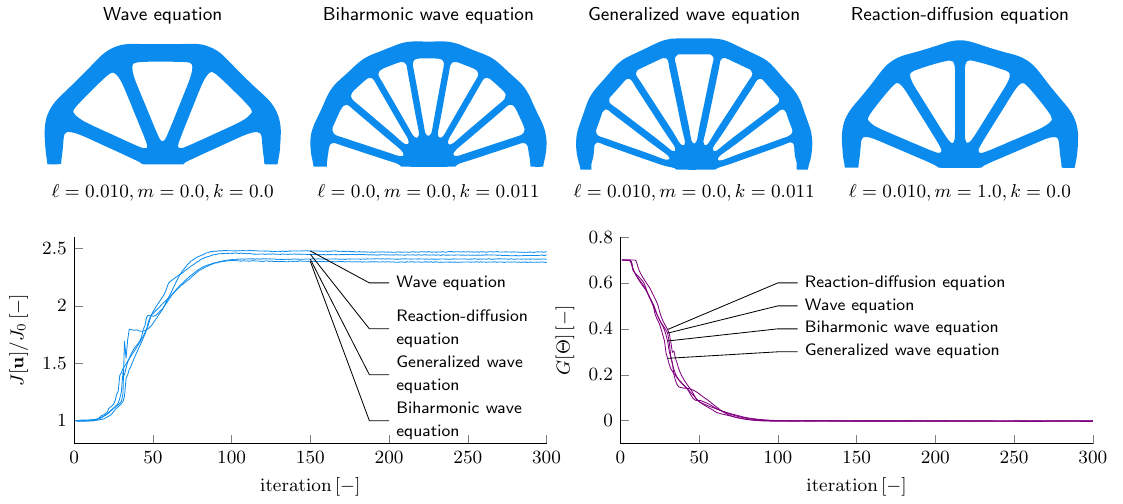}
    \caption{Obtained structural layout using different physical equations and according convergence histories.}
    \label{fig:girder}
\end{figure}

Finally, the methods under consideration are compared in terms of their resulting topologies and structural characteristics using the example of a simply supported girder illustrated in Fig.~(\ref{fig:testcases}c). The parameters are fixed with $\ell = 0.010$, $m = 1.0$ and $k = 0.011$, whereby the respective equations of the evolution problem are clearly determined. Both the standard wave equation and the reaction-diffusion equation lead to moderate topologies. However, the latter leads to a slightly higher complexity in the form of additional holes, as the parameter $\ell$ shows in this case greater sensitivity to the nucleation of new holes compared to the standard wave equation. In contrast, the biharmonic and the generalized wave equation generate significantly more complex structures with almost twice as many holes. Especially in the generalized wave equation, the effect of strut tapering occurs again, which can be observed by bulbous struts in the middle and tapered geometries at the support points, see Fig.~(\ref{fig:girder}).
With regard to the convergence histories, all methods show a similar trend. In this example, the biharmonic wave equation provides the best result in terms of the achieved final ratio $J[\mathbf{u}]/J_0$.
\subsection{Compliant mechanism design considering a fixed design domain}
\begin{figure}[ht]
    \centering
    \includegraphics[width=1.\textwidth]{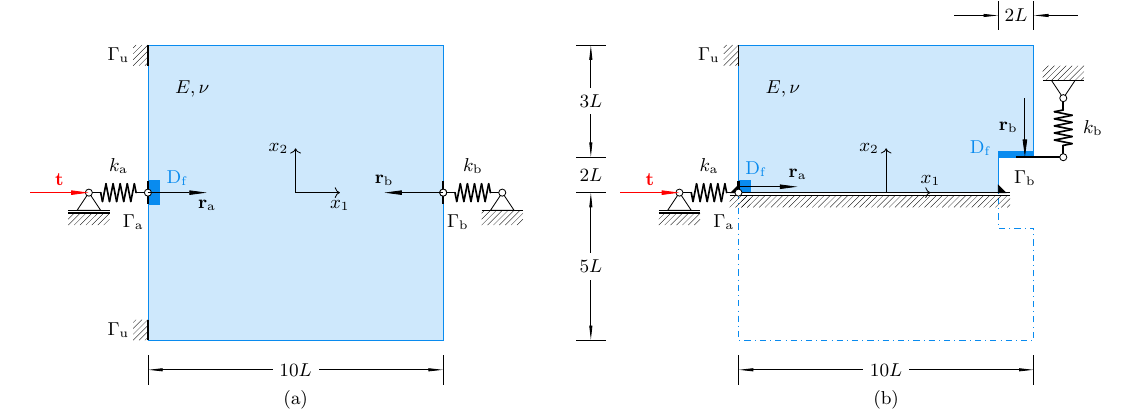}
    \caption{Compliant mechanism design with fixed solid domain $\mathrm{D}_\mathrm{f}$; (a) inverter problem, (b) gripper problem.}
    \label{fig:compliantmechanisms}
\end{figure}

\begin{figure}[ht]
    \centering
    \includegraphics[width=1.\textwidth]{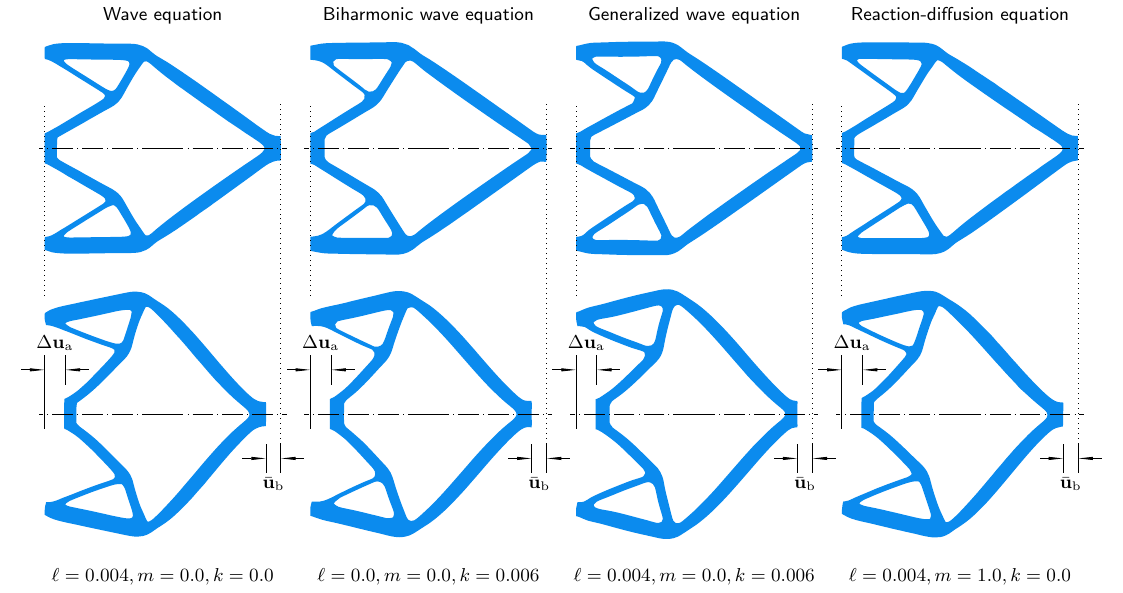}
    \caption{Optimized structural layouts for the inverter problem with respect to different physical equations.}
    \label{fig:inverterresults}
\end{figure}

\begin{figure}[ht]
    \centering
    \includegraphics[width=1.\textwidth]{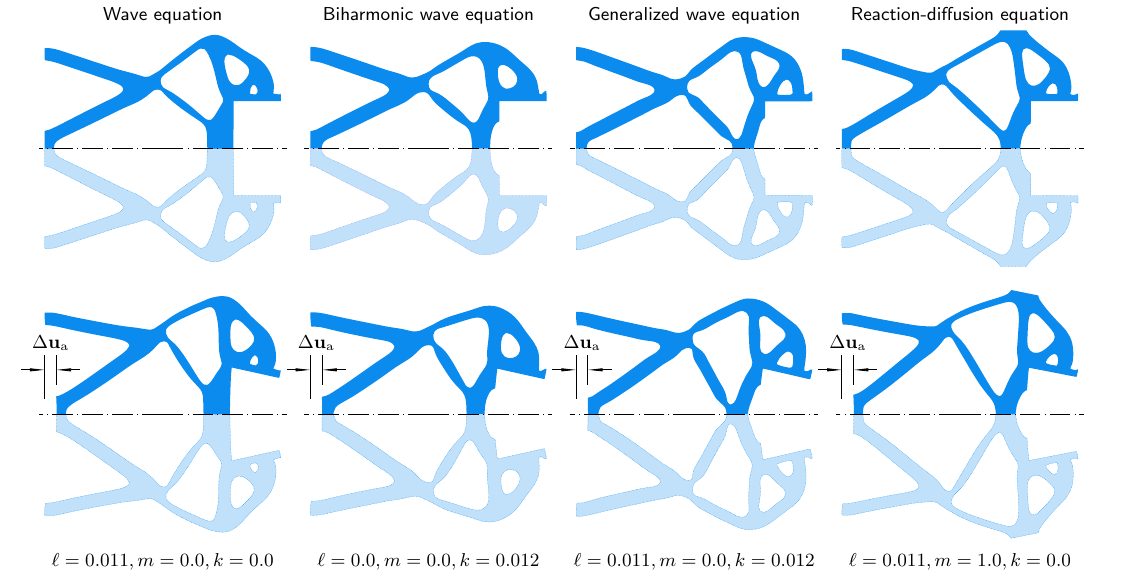}
    \caption{Optimized structural layouts for the gripper problem with respect to different physical equations.}
    \label{fig:gripperresults}
\end{figure}

Compliant mechanisms are elastic structures which, as a result of a defined input displacement or input force, enable desired deformation or force transmission on a given surface and are frequently investigated in structural optimization research \citep{OTORMORI2011,YAMADA2017,EMMENDOERFER2020,EMMENDOERFER2022}. For this purpose, consider a material domain $\Omega$ embedded in a design domain $\mathrm{D}$ containing a predefined domain of solid material $\mathrm{D}_\mathrm{f}$ and being fixed at $\Gamma_\mathrm{u}$. Boundary segment $\Gamma_\mathrm{a}$ denotes the location of the input force application, whereas $\Gamma_\mathrm{b}$ refers to the region of the output deformation. Following \citep{SIGMUND1997}, both boundary segments are equipped with artificial springs of stiffnesses $k_\mathrm{a},k_\mathrm{b}$ that generate reaction forces, thereby facilitating stiffness control. In addition, direction vectors $\mathbf{r}_\mathrm{a},\mathbf{r}_\mathrm{b}\in\mathbb{R}^{N\times 1}$ are incorporated to specify the required direction of input and output deformations. The corresponding optimization problem, which aims to maximize the output deformation under a volume constraint, is then formulated as follows:
\begin{customopti}|s|
    {inf}{\phi\in H^2(\mathscr{D})}{J[\mathbf{u}] = -\int_{\Gamma_\mathrm{b}}\mathbf{r}_\mathrm{b}\cdot\mathbf{u}\,\mathrm{d}\Gamma}{}{}{}
    \addConstraint{G[\Theta]=\frac{1}{V_0}\int_{\mathrm{D}\backslash\mathrm{D}_\mathrm{f}}\Theta\,\mathrm{d}\Omega - V_\mathrm{f}}{\leq 0}{\quad\text{in}\quad\mathrm{D}\backslash\mathrm{D}_\mathrm{f}}{}
    \addConstraint{-\mathrm{div}(\mathbf{C}:\varepsilon(\mathbf{u}))}{=0}{\quad\text{in}\quad\Omega\cup\Omega_\mathrm{f}}{}
    \addConstraint{\mathbf{u}}{=\mathbf{0}}{\quad\text{on}\quad\Gamma_\mathrm{u}}{}
    \addConstraint{(\mathbf{C}:\varepsilon(\mathbf{u}))\mathbf{n}}{=\mathbf{t} - k_\mathrm{a}(\mathbf{r}_\mathrm{a}\otimes\mathbf{r}_\mathrm{a})\mathbf{u}}{\quad\text{on}\quad\Gamma_\mathrm{a}}{}
    \addConstraint{(\mathbf{C}:\varepsilon(\mathbf{u}))\mathbf{n}}{=-k_\mathrm{b}(\mathbf{r}_\mathrm{b}\otimes\mathbf{r}_\mathrm{b})\mathbf{u}}{\quad\text{on}\quad\Gamma_\mathrm{b}}{}
    \addConstraint{(\mathbf{C}:\varepsilon(\mathbf{u}))\mathbf{n}}{=\mathbf{0}}{\quad\text{on}\quad\partial(\Omega\cup\Omega_\mathrm{f})\backslash(\Gamma_\mathrm{u}\cup\Gamma_\mathrm{a}\cup\Gamma_\mathrm{b})}{}.
\end{customopti}
From the boundary conditions, we derive the weak formulation of the state of equilibrium as
\begin{align*}
    R_2:\quad\int_{\mathrm{D}\backslash\mathrm{D}_\mathrm{f}}\varepsilon(\mathbf{u}):\tau(\Theta)\mathbf{C}:\varepsilon(\mathbf{v})\,\mathrm{d}\Omega + \int_{\mathrm{D}_\mathrm{f}}\varepsilon(\mathbf{u}):\mathbf{C}:\varepsilon(\mathbf{v})\,\mathrm{d}\Omega + \dots \\
    \dots + \int_{\Gamma_\mathrm{b}}k_\mathrm{b}(\mathbf{r}_\mathrm{b}\otimes\mathbf{r}_\mathrm{b})\mathbf{u}\cdot\mathbf{v}\,\mathrm{d}\Gamma = \int_{\Gamma_\mathrm{a}}(\mathbf{t} - k_\mathrm{a}(\mathbf{r}_\mathrm{a}\otimes\mathbf{r}_\mathrm{a})\mathbf{u})\cdot\mathbf{v}\,\mathrm{d}\Gamma.   
\end{align*}
while the adjoint equation follows to
\begin{align*}
    \mathcal{A}_2:\quad -\int_{\Gamma_\mathrm{b}}\mathbf{r}_\mathrm{b}\cdot\delta\mathbf{u}\,\mathrm{d}\Gamma = \int_{\mathrm{D}\backslash\mathrm{D}_\mathrm{f}}\varepsilon(\delta\mathbf{u}):\tau(\Theta)\mathbf{C}:\varepsilon(\mathbf{v})\,\mathrm{d}\Omega + \int_{\mathrm{D}_\mathrm{f}}\varepsilon(\delta\mathbf{u}):\mathbf{C}:\varepsilon(\mathbf{v})\,\mathrm{d}\Omega + \dots \\
    \dots + \int_{\Gamma_\mathrm{b}}k_\mathrm{b}(\mathbf{r}_\mathrm{b}\otimes\mathbf{r}_\mathrm{b})\delta\mathbf{u}\cdot\mathbf{v}\,\mathrm{d}\Gamma + \int_{\Gamma_\mathrm{a}}k_\mathrm{a}(\mathbf{r}_\mathrm{a}\otimes\mathbf{r}_\mathrm{a})\delta\mathbf{u}\cdot\mathbf{v}\,\mathrm{d}\Gamma.
\end{align*}
The partial functional derivative of the constraint functional in the augmented formulation with respect to $\Theta$ corresponds to the expression provided in Eq.~(\ref{eq:partialderivativevolumeconstraint}) while the normalization factor is obtained from Eq.~(\ref{eq:standardnormalizationfactor}). Then, the pertubation term results to
\begin{align*}
    f_2 = (\lambda + r G[\Theta])\frac{1}{V_0}\mathbb{I}_G - \varepsilon(\mathbf{u}):\frac{\partial\tau}{\partial\Theta}\mathbf{C}:\varepsilon(\mathbf{v}). 
\end{align*}
We now apply our method to the inverter and the gripper problems shown in Fig.~(\ref{fig:compliantmechanisms}). 
Here, we set Young's modulus to $E = 1\,\mathrm{N / mm^2}$, Poisson's ratio to $\nu = 0.3$ and the length to $L=0.1\,\mathrm{mm}$. The direction vector at the input section is set to $\mathbf{r}_\mathrm{a} = [1,0]^\top$ in both cases while the output direction vector is set to $\mathbf{r}_\mathrm{b} = [-1,0]^\top$ for the inverter problem and to $\mathbf{r}_\mathrm{b} = [0,-1]^\top$ for the gripper problem. The required volume fraction is $V_\mathrm{f} = 0.20$ for the inverter and $V_\mathrm{f} = 0.30$ for the gripper problem. The finite element size ranges between $h_\mathrm{min} = 2.5\times 10^{-4}$ and $h_\mathrm{max} = 2.0 \times 10^{-2}$ and the spring stiffnesses are set to $k_\mathrm{a}=1.0\times 10^5\,\mathrm{N}/\mathrm{mm}$ and $k_\mathrm{b}=1.0\times 10^3\,\mathrm{N}/\mathrm{mm}$. The transition phase parameter is set to $\beta=2.0$. In this context, Fig.~(\ref{fig:inverterresults}) shows the different results with respect to the individual physical equations as well as the deformed shapes under a fixed input displacement $\Delta\mathbf{u}_\mathrm{a}$ in horizontal direction.
We observe that the resulting topologies show only minor variations in shape, and the average output deflection calculated by
\begin{align}
    \bar{\mathbf{u}}_\mathrm{b} = \frac{1}{S_\mathrm{b}}\int_{\Gamma_\mathrm{b}}\mathbf{r}_\mathrm{b}\cdot\mathbf{u}\,\mathrm{d}\Gamma\quad\text{with}\quad S_\mathrm{b} = \int_{\Gamma_\mathrm{b}}\,\mathrm{d}\Gamma 
\end{align} 
remains nearly identical across all cases. Here, the ratio between the output deformation and the input displacement $\bar{\mathbf{u}}_\mathrm{b}/\Delta\mathbf{u}_\mathrm{a}$ results to $0.69$ for the wave equation, $0.66$ for the biharmonic wave equation, $0.70$ for the generalized wave equation and $0.68$ for the reaction-diffusion equation. For the gripper design, the symmetry of the problem is conveniently exploited so that only the upper half of the structure is considered. The final results are shown in Fig.~(\ref{fig:gripperresults}). 
Compared to the inverter problem, there are clear differences in the resulting topologies. While the standard, the generalized wave equation, and the reaction-diffusion equation lead to comparable complex topologies, it is also noticeable with the reaction-diffusion equation that the upper part of the gripper is attracted to the upper boundary of the design domain. In contrast, the other methods produce rounded, pincer-like structures. It is also noticeable that with the biharmonic and the generalized wave equation, the strut thicknesses remain relatively constant, while the standard wave equation and the reaction-diffusion equation show stronger variations in the strut geometry. The bulbous shape of the struts can also be seen here when applying the generalized wave equation. With regard to the ratio between output and input displacement $\bar{\mathbf{u}}_\mathrm{b}/\Delta\mathbf{u}_\mathrm{a}$ we find for the wave equation $0.86$, for the biharmonic wave equation $0.94$, for the generalized wave equation $0.97$, and for the reaction-diffusion equation $0.81$. 
\subsection{Bi-objective self-weight and strain energy minimization under local stress constraint}
\begin{figure}[ht]
    \centering
    \includegraphics[width=1.\textwidth]{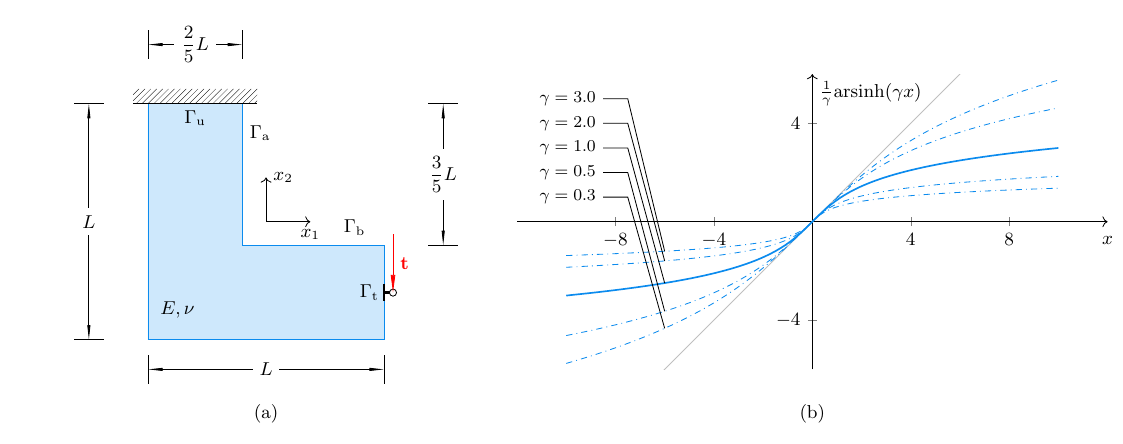}
    \caption{L-bracket problem; (a) problem definition, (b) behavior of inverse hyperbolic sine function with respect to scaling parameter $\gamma$.}
    \label{fig:lbracketproblem}
\end{figure}

\begin{figure}[ht]
    \centering
    \includegraphics[width=1.\textwidth]{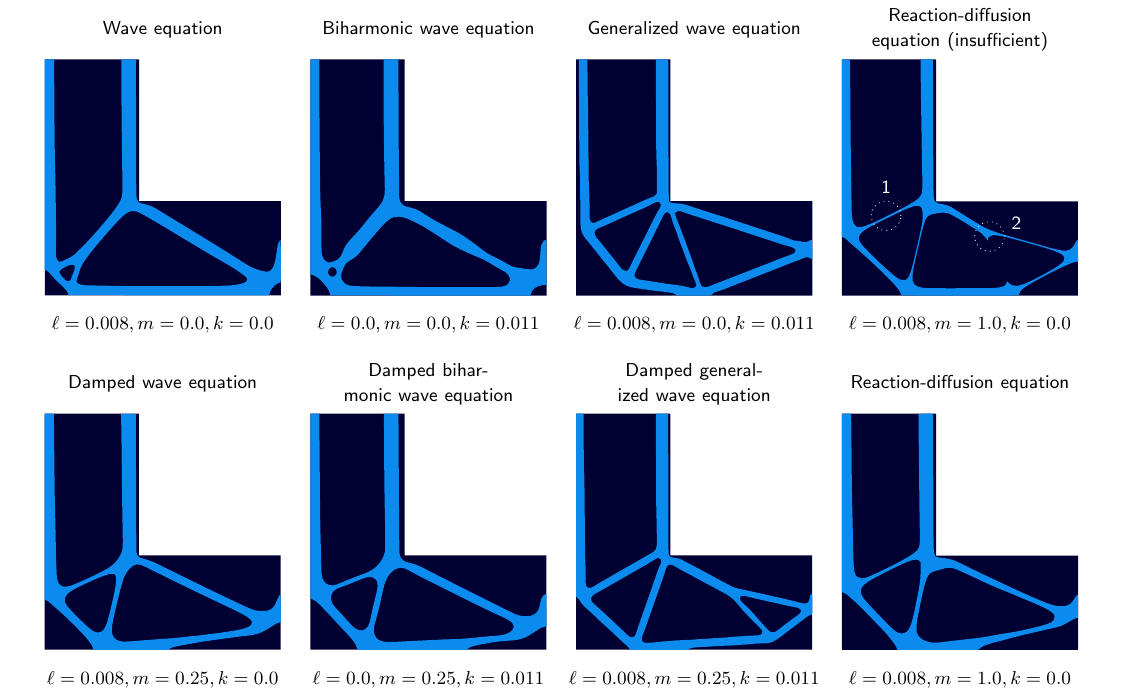}
    \caption{Computed structural layouts with respect to different physical equations.}
    \label{fig:lbrackettopology}
\end{figure}

\begin{figure}[H]
    \centering
    \includegraphics[width=1.\textwidth]{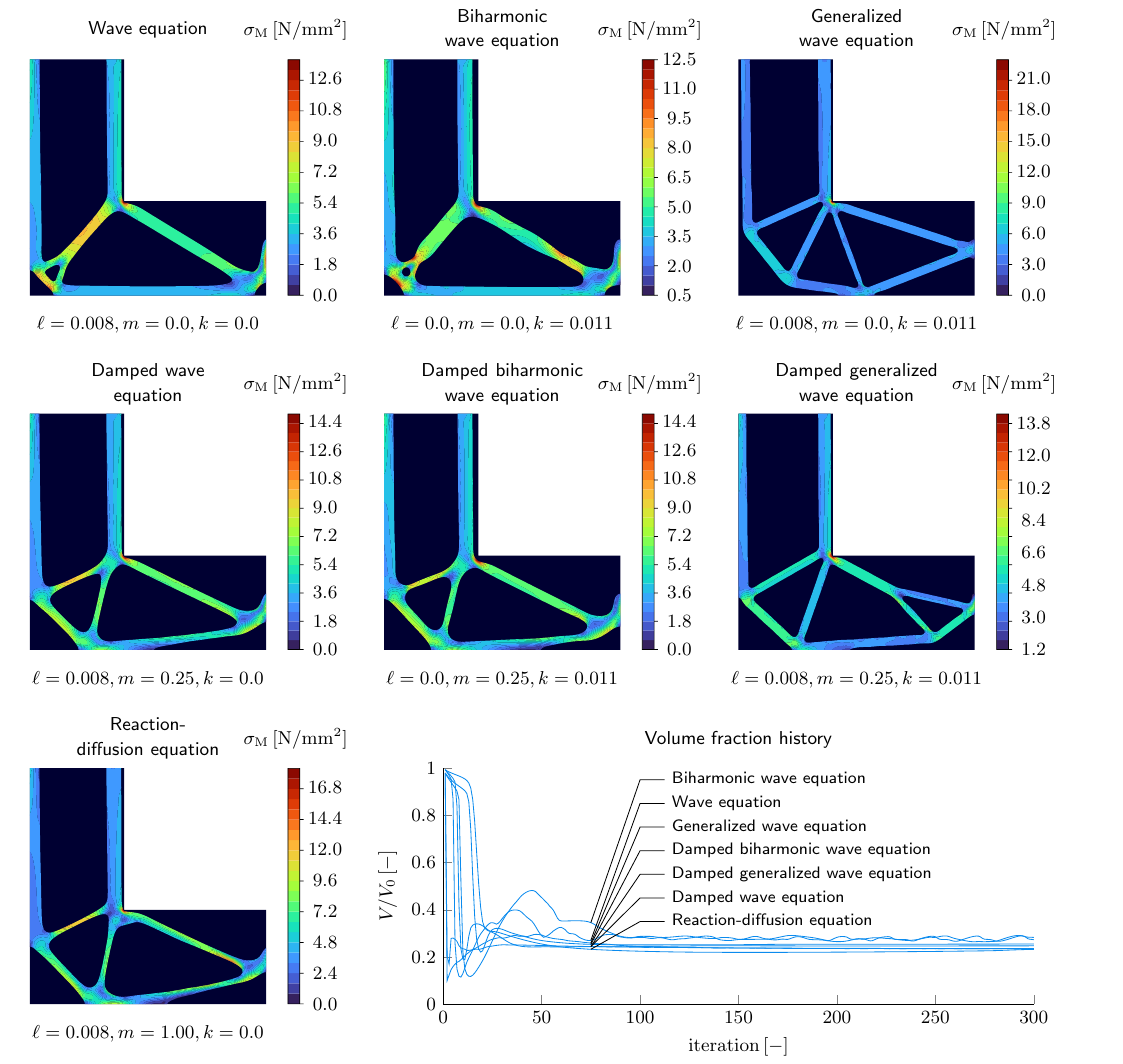}
    \caption{von Mises stress distributions and volume fraction histories with respect to the different physical equations.}
    \label{fig:lbracketstress}
\end{figure}

The consideration of stress constraints in structural optimization is a widely studied problem. A common method is the use of a $p$-norm formulation, in which the stress distribution is aggregated to a global scalar quantity \citep{LE2010,XIA2012,POLAJNAR2017,KAMBAMPATI2021}. This allows the sensitivity to stress peaks to be controlled via the parameter $p$. A major advantage is that the stress information can be represented compactly and integrated efficiently into optimization algorithms. However, local details of the stress distribution are lost due to the aggregation, meaning that the $p$-norm generally cannot guarantee that material failure will not occur at any point in the structure. To address this issue, additional methods are often employed to approximate the maximum stress more accurately \citep{LEE2016,PICELLI2018}. An alternative approach incorporates the consideration of a local stress constraint formulation as shown in \citep{WANG2013,EMMENDOERFER2014,CHENG2024}. In order to demonstrate the broad applicability of our approach we address the problem of bi-objective self-weight and compliance minimization under a local stress constraints. In case of neglected body forces, the corresponding optimization problem is then formulated as follows:
\begin{customopti*}|s|
    {inf}{\phi\in H^2(\mathscr{D})}{J[\mathbf{u},\Theta;w] = \frac{w}{V_0}\int_\mathrm{D}\Theta\,\mathrm{d}\Omega + \frac{1-w}{2}\int_\mathrm{\mathrm{D}}\varepsilon(\mathbf{u}):\tau(\Theta)\mathbf{C}:\varepsilon(\mathbf{u})\,\mathrm{d}\Omega}{}{}{}
    \addConstraint{g(\mathbf{u})=\frac{\sigma_\mathrm{M}(\mathbf{u})}{f_\mathrm{y}} - 1}{\leq 0}{\quad\text{in}\quad\Omega}{}
    \addConstraint{-\mathrm{div}(\mathbf{C}:\varepsilon(\mathbf{u}))}{=\mathbf{0}}{\quad\text{in}\quad\Omega}{}    
    \addConstraint{\mathbf{u}}{=\mathbf{0}}{\quad\text{on}\quad\Gamma_\mathrm{u}}{} 
    \addConstraint{(\mathbf{C}:\varepsilon(\mathbf{u}))\mathbf{n}}{=\mathbf{t}}{\quad\text{on}\quad\Gamma_\mathrm{t}}{} 
    \addConstraint{(\mathbf{C}:\varepsilon(\mathbf{u}))\mathbf{n}}{=\mathbf{0}}{\quad\text{on}\quad\partial\Omega\backslash(\Gamma_\mathrm{u}\cup\Gamma_\mathrm{t})}{}        
\end{customopti*}
with weighting factor $w\in(0,1)$ and $f_\mathrm{y}$ as the maximum permissible stress. In this context, the compliance term ensures the avoidance of structural weak solutions. The von Mises stress is defined by
\begin{align*}
    \sigma_\mathrm{M}(\mathbf{u}) = \sqrt{\frac{3}{2}s(\mathbf{u}):s(\mathbf{u})}
\end{align*}
where $s(\mathbf{u})$ denotes the stress deviator tensor
\begin{align*}
    s(\mathbf{u}) = (\mathbf{C}:\varepsilon(\mathbf{u})) - \frac{1}{3}\mathrm{tr}(\mathbf{C}:\varepsilon(\mathbf{u}))\mathbf{I}
\end{align*}
with $\mathbf{I}$ being the identity tensor of second rank. From the boundary conditions we identify the governing equation 
\begin{align*}
    R_3:\quad\int_\mathrm{D}\varepsilon(\mathbf{u}):\tau(\Theta)\mathbf{C}:\varepsilon(\mathbf{v})\,\mathrm{d}\Omega = \int_{\Gamma_\mathrm{t}}\mathbf{t}\cdot\mathbf{v}\,\mathrm{d}\Gamma
\end{align*}
in its weak formulation. In comparison to the previously treated problems, the inequality constraint must be evaluated point-wise in $\Omega$. Therefore, we interpret the Lagrangian multiplier as a spatially varying function $\lambda(\mathbf{x})\in H^1(\mathrm{D})$, allowing the constraint functional in the augmented formulation to be written as
\begin{align*}
    G^\ast[\mathbf{u},\Theta;\lambda] = \int_\mathrm{D}\left(\lambda\bar{g}(\mathbf{u},\Theta) + \frac{r}{2}\bar{g}^2(\mathbf{u},\Theta)\right)\,\mathrm{d}\Omega.
\end{align*}
Here, the augmented constraint function is defined as
\begin{align*}
    \bar{g}(\mathbf{u},\Theta) = \max\left\{g(\mathbf{u},\Theta);-\frac{\lambda^\prime}{r}\right\}\quad\text{in}\quad\mathrm{D}.
\end{align*}
Next, we aim to derive the partial functional derivatives with respect to the displacement field $\mathbf{u}$. To this end, it is convenient to express the constraint function first in its original level set formulation by means of the characteristic function $\chi$, thereby allowing us to exploit its specific properties in the further process. The constraint function is then defined as
\begin{align*}
    g(\mathbf{u},\chi) = \frac{1}{f_\mathrm{y}}\sqrt{\frac{3}{2}s(\mathbf{u})\chi:s(\mathbf{u})\chi} - 1 = \frac{\sigma_\mathrm{M}(\mathbf{u})}{f_\mathrm{y}}\chi - 1\quad\text{in}\quad\mathrm{D}.
\end{align*}
The partial functional derivative with respect to $\mathbf{u}$ follows to
\begin{align*}
    \left\langle\frac{\partial G^\ast}{\partial\mathbf{u}},\delta\mathbf{u}\right\rangle &= \frac{\mathrm{d}}{\mathrm{d}\epsilon}\left[\int_\mathrm{D}\left(\lambda\bar{g}(\mathbf{u} + \epsilon\delta\mathbf{u},\chi) + \frac{r}{2}\bar{g}^2(\mathbf{u} + \epsilon\delta\mathbf{u},\chi)\right)\,\mathrm{d}\Omega\right]_{\epsilon=0} \\
    &= \int_\mathrm{D}(\lambda + r\bar{g}(\mathbf{u},\chi))\frac{\partial\bar{g}(\mathbf{u},\chi)}{\partial\mathbf{u}}:s(\delta\mathbf{u})\,\mathrm{d}\Omega
\end{align*}
with 
\begin{align*}
    \frac{\partial\bar{g}(\mathbf{u},\chi)}{\partial\mathbf{u}}:s(\delta\mathbf{u}) = 
    \left\{
    \begin{array}{ll}
        \displaystyle\frac{3\chi}{2f_\mathrm{y}\sigma_\mathrm{M}(\mathbf{u})}(s(\mathbf{u}):s(\delta\mathbf{u})) &\displaystyle\quad\text{if}\quad g(\mathbf{u},\chi) > -\frac{\lambda^\prime}{r}\\[10 pt]
        0 &\quad\text{otherwise}
    \end{array}\right..
\end{align*}
At this point, we again introduce an indicator function defined as
\begin{align*}
    \mathbb{I}_g \vcentcolon = 
    \left\{
    \begin{array}{ll}
        1 &\quad\text{if}\quad\lambda^\prime + r g(\mathbf{u},\chi) > 0 \\[10 pt]
        0 &\quad\text{otherwise}
    \end{array}
    \right.
\end{align*}
which results in the alternative expressions
\begin{align*}
    \bar{g}(\mathbf{u},\chi) &= g(\mathbf{u},\chi)\mathbb{I}_g + \left(-\frac{\lambda^\prime}{r}\right)(1 - \mathbb{I}_g),\\
    \frac{\partial\bar{g}(\mathbf{u},\chi)}{\partial\mathbf{u}}:s(\delta\mathbf{u}) &= \frac{3\chi}{2f_\mathrm{y}\sigma_\mathrm{M}(\mathbf{u})}(s(\mathbf{u}):s(\delta\mathbf{u}))\mathbb{I}_g.
\end{align*}
Substituting of these expressions yields the partial functional derivative of the constraint functional with respect to the displacement field:
\begin{align*}
    \left\langle\frac{\partial G^\ast}{\partial\mathbf{u}},\delta\mathbf{u}\right\rangle = \int_\mathrm{D}\left(\lambda + r\left(\frac{\sigma_\mathrm{M}(\mathbf{u})}{f_\mathrm{y}}\chi - 1\right)\mathbb{I}_g - \lambda^\prime(1 - \mathbb{I}_g)\right)\frac{3\chi}{2f_\mathrm{y}\sigma_\mathrm{M}(\mathbf{u})}(s(\mathbf{u}):s(\delta\mathbf{u}))\mathbb{I}_g\,\mathrm{d}\Omega.
\end{align*}
In this context, we benefit from the properties of the characteristic function and the indicator function by utilizing the identities $\chi^2\equiv\chi$ and $\mathbb{I}_g^2\equiv\mathbb{I}_g$. Then, the partial functional derivative simplifies to
\begin{align*}
    \left\langle\frac{\partial G^\ast}{\partial\mathbf{u}},\delta\mathbf{u}\right\rangle = \int_\mathrm{D}(\lambda + r g(\mathbf{u}))\frac{3\tau(\Theta)}{2f_\mathrm{y}\sigma_\mathrm{M}(\mathbf{u})}(s(\mathbf{u}):s(\delta\mathbf{u}))\mathbb{I}_g\,\mathrm{d}\Omega
\end{align*}
where we conveniently replace the characteristic function by the auxiliary function $\tau(\Theta)$. The indicator function is then accordingly expressed as
\begin{align}
    \mathbb{I}_g\vcentcolon =
    \left\{
    \begin{array}{ll}
        1 &\quad\text{if}\quad\lambda^\prime + r g(\mathbf{u},\Theta) > 0\\[10 pt]
        0 &\quad\text{otherwise}
    \end{array}\right.
    \label{eq:indicatorfunctionstress}
\end{align}
while the augmented multiplier function is updated by
\begin{align*}
    \lambda = \max\{\lambda^\prime + rg(\mathbf{u},\Theta);0\}.
\end{align*}
Next, we calculate the partial functional derivative of the objective functional with respect to $\mathbf{u}$:
\begin{align*}
    \left\langle\frac{\partial J}{\partial\mathbf{u}},\delta\mathbf{u}\right\rangle &= \frac{\mathrm{d}}{\mathrm{d}\epsilon}\left[\frac{w}{V_0}\int_\mathrm{D}\Theta\,\mathrm{d}\Omega + \frac{1 - w}{2}\int_\mathrm{D}\varepsilon(\mathbf{u} + \epsilon\delta\mathbf{u}):\tau(\Theta)\mathbf{C}:\varepsilon(\mathbf{u} + \epsilon\delta\mathbf{u})\,\mathrm{d}\Omega\right]_{\epsilon=0} \\
    &= (1-w)\int_\mathrm{D}\varepsilon(\delta\mathbf{u}):\tau(\Theta)\mathbf{C}:\varepsilon(\mathbf{u})\,\mathrm{d}\Omega.
\end{align*}
The derivative of the weak form of the equilibrium condition follows to
\begin{align*}
    \left\langle\frac{\partial R_3}{\partial\mathbf{u}},\delta\mathbf{u}\right\rangle &= \frac{\mathrm{d}}{\mathrm{d}\epsilon}\left[\int_\mathrm{D}\varepsilon(\mathbf{u} + \epsilon\delta\mathbf{u}):\tau(\Theta)\mathbf{C}:\varepsilon(\mathbf{v})\,\mathrm{d}\Omega - \int_{\Gamma_\mathrm{t}}\mathbf{t}\cdot\mathbf{v}\,\mathrm{d}\Gamma\right]_{\epsilon=0} \\
    &= \int_\mathrm{D}\varepsilon(\delta\mathbf{u}):\tau(\Theta)\mathbf{C}:\varepsilon(\mathbf{v})\,\mathrm{d}\Omega.
\end{align*}
The calculated partial functional derivatives can now be used to derive the adjoint equation given as
\begin{align*}
    \mathcal{A}_3:\quad\int_\mathrm{D}(\lambda + r g(\mathbf{u}))\frac{3\tau(\Theta)}{2f_\mathrm{y}\sigma_\mathrm{M}(\mathbf{u})}(s(\mathbf{u}):s(\delta\mathbf{u}))\mathbb{I}_g\,\mathrm{d}\Omega + \dots \\
    \dots + (1-w)\int_\mathrm{D}\varepsilon(\delta\mathbf{u}):\tau(\Theta)\mathbf{C}:\varepsilon(\mathbf{u})\,\mathrm{d}\Omega = \int_\mathrm{D}\varepsilon(\delta\mathbf{u}):\tau(\Theta)\mathbf{C}:\varepsilon(\mathbf{v})\,\mathrm{d}\Omega.
\end{align*}
In the following step, we determine the partial functional derivatives with respect to $\Theta$. For the constraint functional we find
\begin{align*}
    \left\langle\frac{\partial G^\ast}{\partial\Theta},\delta\Theta\right\rangle &= \frac{\mathrm{d}}{\mathrm{d}\epsilon}\left[\int_\mathrm{D}\left(\lambda\bar{g}(\mathbf{u},\Theta + \epsilon\delta\Theta) + \frac{r}{2}\bar{g}^2(\mathbf{u},\Theta + \epsilon\delta\Theta)\right)\,\mathrm{d}\Omega\right]_{\epsilon=0} \\
    &= \int_\mathrm{D}(\lambda + r\bar{g}(\mathbf{u},\Theta))\frac{\partial\bar{g}(\mathbf{u},\Theta)}{\partial\Theta}\delta\Theta\,\mathrm{d}\Omega
\end{align*}
with 
\begin{align*}
    \frac{\partial\bar{g}(\mathbf{u},\Theta)}{\partial\Theta} =
    \left\{
    \begin{array}{ll}
        \displaystyle\frac{\partial\tau}{\partial\Theta}\frac{\sigma_\mathrm{M}(\mathbf{u})}{f_\mathrm{y}} &\displaystyle\quad\text{if}\quad g(\mathbf{u},\Theta) > -\frac{\lambda^\prime}{r} \\[10 pt] 
        0 &\quad\text{otherwise}
    \end{array}\right..
\end{align*}
Here, we again utilize the indicator function given in Eq.~(\ref{eq:indicatorfunctionstress}) and obtain
\begin{align*}
    \frac{\partial\bar{g}(\mathbf{u},\Theta)}{\partial\Theta} = \frac{\partial\tau}{\partial\Theta}\frac{\sigma_\mathrm{M}(\mathbf{u})}{f_\mathrm{y}}\mathbb{I}_g
\end{align*}
which yields the expression
\begin{align*}
    \left\langle\frac{\partial G^\ast}{\partial\Theta},\delta\Theta\right\rangle = \int_\mathrm{D}\left(\lambda + r\left(\frac{\sigma_\mathrm{M}(\mathbf{u})}{f_\mathrm{y}}\tau(\Theta) - 1\right)\mathbb{I}_g - \lambda^\prime(1 - \mathbb{I}_g)\right)\frac{\partial\tau}{\partial\Theta}\frac{\sigma_\mathrm{M}(\mathbf{u})}{f_\mathrm{y}}\mathbb{I}_g\delta\Theta\,\mathrm{d}\Omega.
\end{align*}
Since this is valid for all test functions $\delta\Theta$, the derivative follows to
\begin{align*}
    \frac{\partial G^\ast}{\partial\Theta} = (\lambda + r g(\mathbf{u},\Theta))\frac{\partial\tau}{\partial\Theta}\frac{\sigma_\mathrm{M}(\mathbf{u})}{f_\mathrm{y}}\mathbb{I}_g.
\end{align*}
The partial functional derivative of the objective functional with respect to $\Theta$ follows to
\begin{align*}
    \left\langle\frac{\partial J}{\partial\Theta},\delta\Theta\right\rangle &= \frac{\mathrm{d}}{\mathrm{d}\epsilon}\left[\frac{w}{V_0}\int_\mathrm{D}(\Theta + \epsilon\delta\Theta)\,\mathrm{d}\Omega + \frac{1 - w}{2}\int_\mathrm{D}\varepsilon(\mathbf{u}):\tau(\Theta + \epsilon\delta\mathbf{u})\mathbf{C}:\varepsilon(\mathbf{u})\,\mathrm{d}\Omega\right]_{\epsilon=0} \\
    &= \frac{w}{V_0}\int_\mathrm{D}\delta\Theta\,\mathrm{d}\Omega + \frac{1 - w}{2}\int_\mathrm{D}\varepsilon(\mathbf{u}):\frac{\partial\tau}{\partial\Theta}\mathbf{C}:\varepsilon(\mathbf{u})\delta\Theta\,\mathrm{d}\Omega
\end{align*}
where we identify the corresponding derivative 
\begin{align*}
    \frac{\partial J}{\partial\Theta} = \frac{w}{V_0}+ \frac{1 - w}{2}\varepsilon(\mathbf{u}):\frac{\partial\tau}{\partial\Theta}\mathbf{C}:\varepsilon(\mathbf{u}). 
\end{align*}
The partial functional derivative of the equilibrium condition with respect to $\Theta$ is given as
\begin{align*}
    \left\langle\frac{\partial R_3}{\partial\Theta},\delta\Theta\right\rangle &= \frac{\mathrm{d}}{\mathrm{d}\epsilon}\left[\int_\mathrm{D}\varepsilon(\mathbf{u}):\tau(\Theta + \epsilon\delta\Theta)\mathbf{C}:\varepsilon(\mathbf{v})\,\mathrm{d}\Omega - \int_{\Gamma_\mathrm{t}}\mathbf{t}\cdot\mathbf{v}\,\mathrm{d}\Gamma\right]_{\epsilon=0} \\
    &= \int_\mathrm{D}\varepsilon(\mathbf{u}):\frac{\partial\tau}{\partial\Theta}\mathbf{C}:\varepsilon(\mathbf{v})\delta\Theta\,\mathrm{d}\Omega
\end{align*}
yielding the derivative
\begin{align*}
    \frac{\partial R_3}{\partial\Theta} = \varepsilon(\mathbf{u}):\frac{\partial\tau}{\partial\Theta}\mathbf{C}:\varepsilon(\mathbf{v}).
\end{align*}
Finally, the pertubation term can be formulated as follows:
\begin{align*}
    f_3 = \left\{(\lambda + r g(\mathbf{u},\Theta))\frac{\sigma_\mathrm{M}(\mathbf{u})}{f_\mathrm{y}}\mathbb{I}_g+\varepsilon(\mathbf{u}):\mathbf{C}:\left(\frac{1-w}{2}\varepsilon(\mathbf{u})-\varepsilon(\mathbf{v})\right)\right\}\frac{\partial\tau}{\partial\Theta} + \frac{w}{V_0}.
\end{align*}
In the context of local stress constraints, it is well known that the nucleation of new holes or the separation of material phases leads to pronounced discontinuities in the stress distribution. These discontinuities further affect the perturbation term and may cause the phenomenon of boundary locking, where only certain parts of the boundary evolve while others remain fixed. Inspired by Emmendoerfer and Francello \citep{EMMENDOERFER2016}, we first apply a nonlinear scaling of the pertubation term using the inverse hyperbolic sine function. The parameter $\gamma$ can be used to control the scaling intensity, see Fig.~(\ref{fig:lbracketproblem}b).
Furthermore, the scaled pertubation term is regularized by means of a Helmholtz type partial differential equation \citep{KAWAMOTO2011}, i.e.
\begin{align*}
    \left\{
    \begin{array}{rll}
        -\kappa\nabla^2\bar{f}_3 + \bar{f}_3 &= f_3 &\quad\text{in}\quad\mathrm{D} \\[10 pt]
        \nabla\bar{f}_3\cdot\mathbf{n}&=0 &\quad\text{on}\quad\partial\mathrm{D}
    \end{array}\right.,
\end{align*}
where $\bar{f}_3$ denotes the regularized pertubation term and $\kappa$ refers to the regularization parameter. The weak form is obtained by multiplying with an appropriate test function $\delta\bar{f}\in H^1(\mathrm{D})$, leading to
\begin{align*}
    \int_\mathrm{D}(\kappa\nabla\bar{f}_3\cdot\nabla\delta\bar{f} + \bar{f}_3\delta\bar{f})\,\mathrm{d}\Omega = \int_\mathrm{D}\frac{1}{\gamma}\mathrm{arsinh}(\gamma f_3)\delta\bar{f}\,\mathrm{d}\Omega.
\end{align*}
In this example, the normalization factor is defined as
\begin{align*}
    C_v^\ast = \frac{1}{\mathrm{c}_f V_0}\int_\mathrm{D}\left\vert\bar{f}_3\right\vert_{L^2(\mathrm{D})}\,\mathrm{d}\Omega.
\end{align*}
Consider the L-bracket benchmark problem, shown in Fig.~(\ref{fig:lbracketproblem}a), under the assumption of plane stress conditions. The linear elastic material is characterized by Young's modulus of $E=1\,\mathrm{N}/\mathrm{mm}^2$ and Poisson's ration of $\nu=0.3$. The maximum admissible stress is defined by $f_\mathrm{y}=42\,\mathrm{N}/\mathrm{mm}^2$ and the length is set to $L=1.5\,\mathrm{mm}$. A load magnitude $|\mathbf{t}|=1.0\,\mathrm{N}$ is applied. To account for the presence of adjacent void regions, additional Dirichlet boundary conditions $\phi=-1$ are imposed along $\Gamma_\mathrm{a}$ and $\Gamma_\mathrm{b}$. The finite element size ranges from $h_\mathrm{min} = 1.0\times 10^{-3}$ to $h_\mathrm{max} = 1.0\times 10^{-2}$ and the regularized parameter is set to $\kappa = 1.0\times 10^{-5}$. The transition phase parameter is fixed at $\beta = 5.0$ and the weighting factor is set to $w =0.99$ aiming to prioritize the self-weight minimization. Further numerical parameters can be taken from Tab.~(\ref{tab:stresssettings}).
\begin{table}[ht]
    \centering
    \caption{Numerical settings with respect to the used physical equations.}
    \resizebox{\textwidth}{!}{
    \begin{tabular}{ccccccccc}\toprule
        \textsf{Parameter} & \textsf{WE} & \textsf{DWE} & \textsf{BWE} & \textsf{DBWE} & \textsf{GWE} & \textsf{DGWE} & \textsf{RDE} & \textsf{RDE (insufficient)}\\\midrule
        $\mathrm{c}_f$ & 0.020 & 0.010 & 0.004 & 0.004 & 0.100 & 0.100 & 0.010 & 0.010\\
        $\gamma$ & 4.0 & 4.0 & 3.0 & 3.0 & 3.0 & 3.0 & 4.0 & 3.0 \\\bottomrule
    \end{tabular}
    \label{tab:stresssettings}
    }
\end{table}
Figure~(\ref{fig:lbrackettopology}) shows the optimized topologies obtained with regard to the different physical equations. For both the standard and the biharmonic wave equation, the optimized structures exhibit a pronounced tendency to accumulate material along the lower boundary of the design domain. 
This phenomenon is mitigated by the introduction of a damping term. The application of the generalized wave equation results in the formation of struts with nearly uniform thickness, indicating a consistent regularization effect across the domain. In contrast, the other methods lead to significantly higher variability in strut thicknesses. For the reaction-diffusion equation, a structurally insufficient material distribution results for a scaling parameter $\gamma = 3.0$.  Specifically, the optimized topology contains overly slender struts (see encircled region "1" in Fig.~(\ref{fig:lbracketproblem})) as well as residual material connections that have not been fully removed (encircled region "2"). This behavior is characteristic of reaction-diffusion systems, which were originally developed to model localized crystallization phenomena \citep{CAHN1958,ALLEN1979}. Due to their inherently local dynamics, such models allow only limited spatial propagation of information, thereby impeding global coordination within the design domain. In contrast, wave-based equations enable long-range information transfer \citep{ZACHMANOGLOU1986}, which enables a more balanced and globally coherent material distribution throughout the design domain. Figure~(\ref{fig:lbracketstress}) additionally presents the resulting von Mises stress distributions obtained from the various methods available. 
Among these, the damped biharmonic wave equation proves to be the most effective in smoothing the stress singularity at the critical inner corner of the L-bracket. In contrast, the other approaches result in a comparatively limited and more localized smoothing effect. A further noteworthy observation is the significantly more uniform stress distribution achieved by both the generalized wave equation and its damped version. These methods promote an efficient and balanced stress distribution throughout the structure. Conversely, the remaining methods exhibit markedly heterogeneous stress patterns. Moreover, the convergence histories in Fig.~(\ref{fig:lbracketstress}) clearly demonstrates the beneficial impact of damping in terms of numerical stability. This is reflected in a reduction of oscillations, particularly evident in the results obtained using the standard and the biharmonic wave equations. 
\begin{table}[ht]
    \centering
    \caption{Performance results with respect to the physical equations. $J_1[\Theta]$ refers to the objective functional regarding self-weight, $J_2[\mathbf{u},\Theta]$ refers to the mean compliance.}
    \resizebox{1.\textwidth}{!}{
    \begin{tabular}{ccccccccc}\toprule
        \textsf{Quantity} && \textsf{WE} & \textsf{DWE} & \textsf{BWE} & \textsf{DBWE} & \textsf{GWE} & \textsf{DGWE} & \textsf{RDE} \\\midrule
        $J[\mathbf{u},\Theta;w]/J_0$ & $[-]$ & $0.740$ & $0.662$ & $0.741$ & $0.821$ & $0.386$ & $0.4074$ & $1.148$ \\
        $\max_\Omega\{g(\mathbf{u})\}$ & $[-]$ & $-0.675$ & $-0.6534$ & $-0.710$ & $-0.651$ & $-0.455$ & $-0.660$ & $-0.580$\\
        $J_1[\Theta]$ & $[-]$ & $0.273$ & $0.247$ & $0.283$ & $0.255$ & $0.236$ & $0.235$ & $0.232$ \\
        $J_2[\mathbf{u},\Theta]$ & $[\mathrm{N}\,\mathrm{mm}]$ & $93.54$ & $83.16$ & $91.85$ & $113.56$ & $30.54$ & $34.76$ & $183.20$ \\
        $\max_\Omega\{\sigma_\mathrm{M}(\mathbf{u})\}$ & $[\mathrm{N}/\mathrm{mm}^2]$ & $13.64$ & $14.57$ & $12.20$ & $14.65$ & $22.88$ & $14.27$ & $17.65$ \\\bottomrule
    \end{tabular}
    }
    \label{tab:stressresults}
\end{table}
Table~(\ref{tab:stressresults}) shows that none of the solutions exceed the allowable stress limit. The generalized wave equation yields the best results in terms of the achieved aggregated objective functional ratio $J[\mathbf{u},\Theta;w]/J_0$. However, it also exhibits the highest maximum von Mises stress compared to the remaining solutions.

  \section{Conclusion}
\label{sec:conclusion}
In this paper, we proposed a novel framework for systematically deriving an update scheme for the level set function in topology optimization. The central element of the approach was the construction of an auxiliary domain filled with fictitious matter subject to excitation by external conditions. In this connection, the level set function was interpreted as a generalized coordinate describing the deformation of this matter. After assigning kinetic and potential energy, Hamilton's principle was employed to derive the corresponding equation of motion, resulting in a modified wave equation. The formulation was transferred into a numerical scheme via variational methods, and the proposed approach was validated through various design optimization problems, including a detailed analysis of the influence of numerical parameters. The main achievements of our research are summarized as follows:
\begin{enumerate}
    \item By selecting appropriate parameter combinations, different physical equations can be recovered, including the standard and biharmonic wave equations.
    \item The method integrates seamlessly with the current state of research, as the reaction-diffusion equation emerges as a special case by consideration of pure damping.
    \item Parameter $\beta$ controls the width of the transition phase between void and material domain in the approximated Heaviside function and can be used to improve symmetry properties and convergence speed.
    \item In the standard wave equation, parameter $\ell$ enables control over the topological complexity of the structural layout. Introducing damping via the parameter $m$ facilitates the identification of a preferred topology as the evolution is dragged to this solution.
    \item In the biharmonic wave equation, parameter $k$ also allows for tuning the structural complexity. In highly complex topologies, a non-linear pattern in strut thickness emerges. Additional damping yields similar effects as in the standard wave equation.
    \item Combining $\ell$ and $k$ yields the so-called generalized wave equation. We showed that structurally viable solutions are obtained for a ratio $\ell/k < 1$, while a transition to non-viable solutions occurs near $\ell/k \approx 1$. Furthermore, the resulting struts exhibit bulbous mid-sections tapering toward the joints. This  effect becomes more pronounced near the transition region. Damping has a comparable influence here as well.
    \item The mathematical properties of the wave equation promote a more balanced material distribution than the reaction-diffusion approach, as local information is propagated over greater distances. This tends to suppress the formation of thin struts and helps to avoid freezing of local regions. 
    \item Numerical stability can be improved by adjusting the parameter $m$, which effectively suppresses oscillations, as demonstrated in the problem of bi-objective self-weight and compliance minimization under consideration of a local stress constraint.
\end{enumerate}

  % ... 
  \section*{Acknowlegdement}
  This work was supported by the JSPS KAKENHI Grant Number 24KF0144. 

  % ... 
  \section*{Declaration of generative AI and AI-assisted technologies in the writing process}
  During the preparation of this work the authors used DeepL and ChatGPT for English proofreading. After using this tool, the authors reviewed and edited the content as needed and take full responsibility for the content of the published article.

  \appendix
  \section{Appendix}
\label{appendix}

\subsection{Wave equation}
The evolution problem is defined as
\begin{align}
    \left\{\begin{array}{rll}
    \displaystyle \square\phi &= \displaystyle -\frac{f(\mathbf{u},\Theta;\mathbf{v},\lambda)}{c^2}\beta  &\quad\text{in}\quad\mathscr{D}\\[10pt]
    \phi &= 1 &\quad\text{on}\quad\mathscr{G}_\mathrm{t} \\[10pt]
    \nabla\phi\cdot\mathbf{n} &=0 &\quad\text{on}\quad \partial\mathscr{D}\backslash\mathscr{G}_\mathrm{t} \\[10pt]
    \phi\vert_{v=0} &= \phi_0 \\[5pt]
    \displaystyle\left.\frac{\partial\phi}{\partial v}\right\vert_{v=0} &= \displaystyle\left(\frac{\partial\phi}{\partial v}\right)_0
    \end{array}\right..
    \label{eq:standardwaveequation}
\end{align} 
From the discretized equation
\begin{align*}
    \phi_{v+1} - \ell^2\nabla^2\phi_{v+1}= - S_v + 2\phi_v - \phi_{v-1}
\end{align*}
the weak formulation is obtained:
\footnotesize{
\begin{align*}
    \left\{
    \begin{array}{rll}
        \displaystyle\int_{\mathscr{D}}\phi_{v+1}\delta\phi\,\mathrm{d}\mathscr{D} + \int_{\mathscr{D}}\ell^2\nabla\phi_{v+1}\cdot\nabla\delta\phi\,\mathrm{d}\mathscr{D} &= \displaystyle\int_{\mathscr{D}}(-S_v + 2\phi_v - \phi_{v-1})\delta\phi\,\mathrm{d}\mathscr{D} & \quad\text{in}\quad\mathscr{D} \\[10 pt]
        \phi_{v+1} &=0 &\quad\text{on}\quad\mathscr{G}_\mathrm{t}
    \end{array}\right..
\end{align*}
}\normalsize

\subsubsection{Damped wave equation}
The evolution problem is defined as
\begin{align}
    \left\{\begin{array}{rll}
    \displaystyle \square\phi + \frac{d}{c^2}\frac{\partial\phi}{\partial v} &= \displaystyle -\frac{f(\mathbf{u},\Theta;\mathbf{v},\lambda)}{c^2}\beta  &\quad\text{in}\quad\mathscr{D}\\[10pt]
    \phi &= 1 &\quad\text{on}\quad\mathscr{G}_\mathrm{t} \\[10pt]
    \nabla\phi\cdot\mathbf{n} &=0 &\quad\text{on}\quad \partial\mathscr{D}\backslash\mathscr{G}_\mathrm{t} \\[10pt]
    \phi\vert_{v=0} &= \phi_0 \\[5pt]
    \displaystyle\left.\frac{\partial\phi}{\partial v}\right\vert_{v=0} &= \displaystyle\left(\frac{\partial\phi}{\partial v}\right)_0
    \end{array}\right..
\end{align} 
From the discretized equation
\begin{align*}
    (1 + m)\phi_{v+1} - \ell^2\nabla^2\phi_{v+1} = - S_v + (2 + m)\phi_v - \phi_{v-1}
\end{align*}
the weak formulation is obtained:
\footnotesize{
\begin{align*}
    \left\{
    \begin{array}{rll}
        \displaystyle\int_{\mathscr{D}}(1+m)\phi_{v+1}\delta\phi\,\mathrm{d}\mathscr{D} + \int_{\mathscr{D}}\ell^2\nabla\phi_{v+1}\cdot\nabla\delta\phi\,\mathrm{d}\mathscr{D} &= \displaystyle\int_{\mathscr{D}}(-S_v + (2+m)\phi_v - \phi_{v-1})\delta\phi\,\mathrm{d}\mathscr{D} & \quad\text{in}\quad\mathscr{D} \\[10 pt]
        \phi_{v+1} &=0 &\quad\text{on}\quad\mathscr{G}_\mathrm{t}
    \end{array}\right..
\end{align*}
}\normalsize

\subsection{Biharmonic wave equation}
The evolution problem is defined as
\begin{align}
    \left\{\begin{array}{rll}
    \displaystyle \frac{\partial^2\phi}{\partial v^2} - h^4\nabla^4\phi &= \displaystyle f(\mathbf{u},\Theta;\mathbf{v},\lambda)\beta  &\quad\text{in}\quad\mathscr{D}\\[10pt]
    \phi &= 1 &\quad\text{on}\quad\mathscr{G}_\mathrm{t} \\[10pt]
    \nabla\phi\cdot\mathbf{n} &=0 &\quad\text{on}\quad \partial\mathscr{D}\backslash\mathscr{G}_\mathrm{t} \\[10pt]
    \phi\vert_{v=v_0} &= \phi_0 \\[5pt]
    \displaystyle\left.\frac{\partial\phi}{\partial v}\right\vert_{v=v_0} &= \displaystyle\left(\frac{\partial\phi}{\partial v}\right)_0
    \end{array}\right..
    \label{eq:biharmonicwaveequation}
\end{align} 
From the discretized equations
\begin{align*}
    \phi_{v+1} - k^4\nabla^2\omega &= - S_v + 2\phi_v - \phi_{v-1}, \\
    \nabla^2\phi_{v+1} &= \omega
\end{align*}
the weak formulations of the coupled system of equations are obtained:
\begin{align*}
    \left\{
    \begin{array}{rll}
        \displaystyle\int_\mathscr{D}\phi_{v+1}\delta\phi\,\mathrm{d}\mathscr{D} + \int_\mathscr{D}k^4\nabla\omega\cdot\nabla\delta\phi\,\mathrm{d}\mathscr{D} &= \dots &\\[10pt]
        \displaystyle\dots = \int_\mathscr{D}(-S_v + 2\phi_v - \phi_{v-1})\delta\phi\,\mathrm{d}\mathscr{D} & &\quad\text{in}\quad\mathscr{D}\\[10pt]
        \displaystyle\int_\mathscr{D}\omega\delta\omega\,\mathrm{d}\mathscr{D} + \int_\mathscr{D}\nabla\phi_{v+1}\cdot\nabla\delta\omega\,\mathrm{d}\mathscr{D} &= 0 &\quad\text{in}\quad\mathscr{D}\\[10 pt]
        \phi_{v+1} &= 1&\quad\text{on}\quad\mathscr{G}_\mathrm{t}
    \end{array}\right..
\end{align*}

\subsection{Damped biharmonic wave equation}
The evolution problem is defined as
\begin{align}
    \left\{\begin{array}{rll}
    \displaystyle \frac{\partial^2\phi}{\partial v^2} + d\frac{\partial\phi}{\partial v} - h^4\nabla^4\phi &= \displaystyle f(\mathbf{u},\Theta;\mathbf{v},\lambda)\beta  &\quad\text{in}\quad\mathscr{D}\\[10pt]
    \phi &= 1 &\quad\text{on}\quad\mathscr{G}_\mathrm{t} \\[10pt]
    \nabla\phi\cdot\mathbf{n} &=0 &\quad\text{on}\quad \partial\mathscr{D}\backslash\mathscr{G}_\mathrm{t} \\[10pt]
    \phi\vert_{v=v_0} &= \phi_0 \\[5pt]
    \displaystyle\left.\frac{\partial\phi}{\partial v}\right\vert_{v=v_0} &= \displaystyle\left(\frac{\partial\phi}{\partial v}\right)_0
    \end{array}\right..
    \label{eq:dampedbiharmonicwaveequation}
\end{align} 
From the discretized equations
\begin{align*}
    (1+m)\phi_{v+1} - k^4\nabla^2\omega &= - S_v + (2+m)\phi_v - \phi_{v-1}, \\
    \nabla^2\phi_{v+1} &= \omega
\end{align*}
the weak formulations of the coupled system of equations are obtained:
\begin{align*}
    \left\{
    \begin{array}{rll}
        \displaystyle\int_\mathscr{D}(1+m)\phi_{v+1}\delta\phi\,\mathrm{d}\mathscr{D} + \int_\mathscr{D}k^4\nabla\omega\cdot\nabla\delta\phi\,\mathrm{d}\mathscr{D} &= \dots &\\[10pt]
        \displaystyle\dots = \int_\mathscr{D}(-S_v + (2+m)\phi_v - \phi_{v-1})\delta\phi\,\mathrm{d}\mathscr{D} & &\quad\text{in}\quad\mathscr{D}\\[10pt]
        \displaystyle\int_\mathscr{D}\omega\delta\omega\,\mathrm{d}\mathscr{D} + \int_\mathscr{D}\nabla\phi_{v+1}\cdot\nabla\delta\omega\,\mathrm{d}\mathscr{D} &= 0 &\quad\text{in}\quad\mathscr{D}\\[10 pt]
        \phi_{v+1} &= 1&\quad\text{on}\quad\mathscr{G}_\mathrm{t}
    \end{array}\right..
\end{align*}

\subsection{Generalized wave equation}
The evolution problem is defined as
\begin{align}
    \left\{\begin{array}{rll}
        \displaystyle \square\phi - \frac{h^4}{c^2}\nabla^4\phi &= \displaystyle -\frac{f(\mathbf{u},\Theta;\mathbf{v},\lambda)}{c^2}\beta  &\quad\text{in}\quad\mathscr{D}\\[10 pt]
        \phi &= 1 &\quad\text{on}\quad\mathscr{G}_\mathrm{t} \\[10 pt]
        \nabla\phi\cdot\mathbf{n} &=0 &\quad\text{on}\quad \partial\mathscr{D}\setminus\mathscr{G}_\mathrm{t} \\[10 pt]
        \phi\vert_{v=0} &= \phi_0 \\[10 pt]
        \displaystyle\left.\frac{\partial\phi}{\partial v}\right\vert_{v=0} &= \displaystyle\left(\frac{\partial\phi}{\partial v}\right)_0
        \end{array}\right.
        \label{eq:generalizedwaveequation}
\end{align}
From the discretized equations
\begin{align*}
    \phi_{v+1} - \ell^2\omega - k^4\nabla^2\omega = -S_v + 2\phi_v - \phi_v
\end{align*}
the weak formulation is of the coupled system of equations is obtained:
\footnotesize{
\begin{align*}
    \left\{
    \begin{array}{rll}
        \displaystyle\int_\mathscr{D}\phi_{v+1}\delta\phi\,\mathrm{d}\mathscr{D} - \int_\mathscr{D}\ell^2\omega\delta\phi\,\mathrm{d}\mathscr{D} + \int_\mathscr{D}k^4\nabla\omega\cdot\nabla\delta\phi\,\mathrm{d}\mathscr{D} &= \dots &\\[10pt]
        \displaystyle\dots = \int_\mathscr{D}(-S_v + 2\phi_v - \phi_{v-1})\delta\phi\,\mathrm{d}\mathscr{D} & &\quad\text{in}\quad\mathscr{D}\\[10pt]
        \displaystyle\int_\mathscr{D}\omega\delta\omega\,\mathrm{d}\mathscr{D} + \int_\mathscr{D}\nabla\phi_{v+1}\cdot\nabla\delta\omega\,\mathrm{d}\mathscr{D} &= 0 &\quad\text{in}\quad\mathscr{D}\\[10 pt]
        \phi_{v+1} &= 1&\quad\text{on}\quad\mathscr{G}_\mathrm{t}
    \end{array}\right.
\end{align*}
}\normalsize

\subsection{Reaction diffusion equation}
The equation of motion is given as
\begin{align}
    \left\{\begin{array}{rll}
    \displaystyle d\frac{\partial\phi}{\partial v} - c^2\nabla^2\phi_{v+1} &= \displaystyle -f(\mathbf{u},\Theta;\mathbf{v},\lambda)\beta  &\quad\text{in}\quad\mathscr{D}\\[10pt]
    \phi &= 1 &\quad\text{on}\quad\mathscr{G}_\mathrm{t} \\[10pt]
    \nabla\phi\cdot\mathbf{n} &=0 &\quad\text{on}\quad \partial\mathscr{D}\backslash\mathscr{G}_\mathrm{t} \\[10pt]
    \phi\vert_{v=0} &= \phi_0 \\
    \end{array}\right.
    \label{eq:rdeequation}
\end{align}
From the discretized equation
\begin{align*}
    m\phi_{v+1} - \ell^2\nabla^2\phi_{v+1} = -S_v + m\phi_v
\end{align*}
the weak formulation is obtained:
\footnotesize{
\begin{align*}
    \left\{
    \begin{array}{rll}
        \displaystyle\int_{\mathscr{D}}m\phi_{v+1}\delta\phi\,\mathrm{d}\mathscr{D} + \int_{\mathscr{D}}\ell^2\nabla\phi_{v+1}\cdot\nabla\delta\phi\,\mathrm{d}\mathscr{D} &= \displaystyle\int_{\mathscr{D}}(-S_v + m\phi_v)\delta\phi\,\mathrm{d}\mathscr{D} &\quad\text{in}\quad\mathscr{D} \\[10 pt]
        \phi_{v+1} &= 1&\quad\text{on}\quad\mathscr{G}_\mathrm{t}
    \end{array}\right..
\end{align*}
}

  \bibliographystyle{elsarticle-num-names} 
  \bibliography{BIB_literature.bib}





\end{document}